%% file: EEq.tex
\documentclass[aps,preprint]{revtex4}

\usepackage{algorithm}
\usepackage{algpseudocode}

\usepackage{makeidx}
\usepackage{multirow}
\usepackage{array}

\usepackage{amsmath,amssymb,amsthm}
\usepackage{mathrsfs}
\usepackage{amsfonts}
\usepackage[]{xcolor}
\usepackage{graphicx}
\usepackage{epstopdf}
\usepackage{ulem}
\usepackage{hyperref}
\newcommand{\ud}{\mathrm{d}}

\begin{document}

\title{On the estimating equations and objective functions for parameters of exponential power distribution: Application for disorder}

	\author{Mehmet Niyazi \c{C}ankaya$^{a1,a2}$}
	
	\affiliation{$^{a1}$ Faculty of Applied Sciences, Department of International Trading and Finance, U\c{s}ak University,  U\c{s}ak, Turkey}
	
	\affiliation{$^{a2}$ Faculty of Art and Sciences, Department of Statistics, U\c{s}ak University, U\c{s}ak, Turkey}

\begin{abstract}
The efficient modeling for disorder in a phenomena depends on the chosen score and  objective functions. The main parameters in modeling are location, scale and shape. The exponential power distribution known as generalized Gaussian is extensively used in modeling. In real world, the observations are member of different parametric models or disorder in a data set  exists.  In this study, estimating equations for the parameters of exponential power  distribution are derived to have robust and also efficient M-estimators when the data set includes disorder or contamination.  The robustness property of M-estimators for the parameters  is examined.  Fisher information matrices based on the derivative of score functions from $\log$, $\log_q$ and distorted log-likelihoods are proposed by use of Tsallis $q$-entropy  in order to have variances of M-estimators. It is shown that matrices derived by score functions are positive semidefinite if conditions are satisfied.  Information criteria inspired by Akaike and Bayesian are arranged by taking the absolute value of score functions. Fitting performances of score functions from estimating equations and objective functions are tested by applying volume, information criteria and mean absolute error which are essential tools in modeling to assess the fitting competence of the proposed functions.  Applications from simulation and real data sets are carried out to compare the performance of estimating equations and objective functions. It is generally observed that the distorted log-likelihood for the estimations of parameters of exponential power distribution has superior performance than other score and objective functions for the contaminated data sets.

\end{abstract}
\maketitle

\section{Introduction}

Efficient modeling a data set is an important problem in the applied science and depends on the functions which are used to perform modeling. If a data set includes contamination or disorder,  it is difficult to manage the efficiency. The relative entropies  or divergences are used to estimate robustly the parameters in a parametric model $f(x;\boldsymbol \theta)$ \cite{PardoSD}  when  contamination exists in the data set. Recently, Tsallis statistic based on $q$-deformed calculus gains the most important concern to model the data set including contamination or disorder ~\cite{Tsallisbook09,Tsallis88,AbeOka01book}. The statistics for location and scale are inevitable indicators to summarize a data set. We assume that a data set is a member of $f(x;\boldsymbol \theta)$, which is strict assumptation to model and summarize the tendeny (location) and dispersion (scale) of the data set. Getting location and scale is performed by using two procedures which are $f(x;\boldsymbol \theta)$ and estimation method for parameters $\boldsymbol{\theta}$ of $f$.   The most  prominent estimation method is maximum likelihood estimation (MLE) which is equivalent to minimization of relative divergences between data and function \cite{PardoSD}. MLE method produces estimators $\hat{\boldsymbol \theta}$ which have properties such as unbiasedness, efficiency, minimum variance \cite{LehmannCas98}. In MLE method, we can take ordinary logarithm, i.e. $\log$, to proceed a simple analytical expression. The chosen parametric function and taking logarithm in MLE is originally an objective function (M-function) $\rho=-\log(f)$ which gives an advantage to manage the efficiency and so we have efficient M-estimators from M-functions. In this case,  $\log$ can be replaced by  $\log_q$ derived by Tsallis $q$-entropy. Thus, we have maximum $\log_q$ likelihood estimation method which has been extensively used and applied by \cite{FerrariYang10,CanKor18,MLqEBio,FerrariOpere} for the estimations of parameters. The parameter $q$ makes a  slow transition to $\log$. The different values of $q$ in $\log_q$ zooms to tail or central parts of the model. It is possible to consider other types of the generalized and deformed logarithms \cite{GentEntReview1,JizbaArimitsu,Bercher12a,Jizba3,JizKorhybrid}. However, the concavity, analytical simplicity for examining finiteness of score functions and constrained optimization in computation for $\log_q$ are important properties to apply for estimations of parameters \cite{CanKor18}. $\log_q$ is simple and free from integral calculations when we compare with divergences \cite{Basuetal98} which is originally from \cite{Vajda86}.  It can be convenient to propose a log-score function from  $\log$ by use of MLE instead of scanning the model $f(x;\boldsymbol \theta)$ completely by means of $\log_q$. Thus, it is possible to model the disorder or different contaminations into  data sets. As an another M-function based on objective function, a convex combination of underlying and contamination distributions are proposed by \cite{CanKor18} for estimations of parameters in the underlying distribution. 

It is assumed that a data set comes from  an underlying distribution $f(x;\boldsymbol \theta)$. In this case, observations $x_1,x_2,\cdots,x_n$ are a member of $f(x;\boldsymbol \theta)$ and so they are identical. However, identically distributed observations are not generally observed in an empirical application. The underlying distribution $f(x;\boldsymbol \theta)$ can contain some contamination from $f(x;\boldsymbol \tau)$. If there exists a contamination, using robust estimation methods for parameters of underlying distribution is necessary \cite{Hub64}. In order to make a robust and also an efficient fitting for non-identical observations, the estimating equations (EEs) from a generalized version of likelihood estimation methods can be used \cite{God60}. EEs are defined as solutions of system of equations according to corresponding parameters $\boldsymbol \theta$: 

\begin{equation}\label{defofpsiFintro}
\sum_{i=1}^n \Psi(x_i; \boldsymbol \theta)_{|\boldsymbol \theta=\hat{\boldsymbol \theta}}=\boldsymbol 0,
\end{equation}
where $\Psi(x_i;\boldsymbol \theta)=\frac{\partial}{\partial \boldsymbol \theta} \Lambda[f(x_i;\boldsymbol \theta)]$, $\Lambda$ represents a generalization of $\log$ and $\Psi$ is a vector for score functions $\psi_1, \psi_2, \cdots \psi_p$. M-estimators $\hat{\boldsymbol \theta}$ are produced by EEs which include M-functions \cite{CanKor18,God60}. If $\Lambda$ is $\log$, then MLE of parameters are obtained. One can find different examples of M-functions from  \cite{Hub64,Omer98,Hampeletal86}. The parameters of underlying distribution are estimated robustly by use of these M-functions.  However, efficiency, that is, the more precise fitting performed by M-functions, is another challenging task in the estimations of parameters. For this reason, it is necessary to propose M-functions which can accomplish to manage the efficiency. The M-functions from  objective and score can be chosen by volume from information geometry, information criteria (IC) and mean absolute error.  \cite{CanKor18,infocomplestamod,MyungBalasubramanian,ClusteringFisher,Maybankvol}.

This paper aims to propose the score function $S$ which can be generated from $\Psi(x; \boldsymbol \theta)$ by means of EEs.  We consider to derive new score functions $S_q$ and $S^D$ from $\log_q$ and distorted log likelihoods, respectively. Thus, a weighted log-score function $wS$ from $S_q$ and $S^D$ is obtained. Since M-estimators from M-functions are competative each other, we will make a comparison among them to test their performance about getting the efficient and the robust M-estimators.  Fisher information and IC for score functions will be proposed to perform inference based on score functions.

The remainder of this paper is organized as follows. Section \ref{Tsallissection} introduces Tsallis $q$-entropy. Section \ref{Llocscale} introduces maximum likelihood type estimation methods.  In Section \ref{prelimi}, we briefly recall definitions of M-estimation originally generated from MLE. EEs from likelihood estimation are proposed.  Thus, we will propose score functions from EEs for different values $\alpha_j$, $j=1,2,3$ of shape parameter $\alpha$ at some intervals on the real line and so the location and scale parameters of underlying distribution can be estimated efficiently and robustly. Fisher information  are given by Section \ref{FisherEEs}. Section \ref{volumeandIC} provides tools for information theory.  Real data application is given by Section \ref{simrealapp}.  Section \ref{labconclusions}  is divided into conclusion and future works. Appendices are given for random number generation, computation, simulation, positive semidefinite of Fisher information matrices and the corresponding results of these matrices for EP distribution.  

\section{$q$-derivative and Tsallis $q$-entropy}\label{Tsallissection}
Entropy is a tool to measure disorder in a phenomena \cite{AmariIG}. Tsallis $q$-entropy is a generalization of Shannon entropy \cite{Shannon} based on $q$-deformed calculus \cite{Ernst12}.  The idea which generates a functional form of entropy is provided by \cite{Abe97}. There are two steps to produce Tsallis $q$-entropy. Firstly, we have partition function given by
\begin{equation}\label{partitonf}
g(z;\boldsymbol{\theta})=f(x;\boldsymbol{\theta})^z,
\end{equation}
\noindent where  $z \in \mathbb{R}$ and vector $\boldsymbol \theta \in \mathbb{R}^p$ represents the number of $p$ parameters $\theta_1, \theta_2$ and $\theta_p$, $p \in \mathbb{Z}^{+}$ from a parametric model $f(x;\boldsymbol \theta) \in \mathbb{R}$. Secondly, we have Jackson $q$-derivative defined as 
\begin{equation}\label{Jacksonq}
D_z^q g(z;\boldsymbol{\theta})=-\frac{g(z;\boldsymbol{\theta})-g(qz;\boldsymbol{\theta})}{(1-q)z}, ~0<q<1.
\end{equation}
When the definition in Eq. \eqref{Jacksonq} is applied to Eq. \eqref{partitonf}, the functional form of Tsallis $q$-entropy is produced by the following form: 
\begin{equation}\label{1qz}
D_z^q f(x;\boldsymbol{\theta})^z|_{z=1}=-\frac{f(x;\boldsymbol \theta)-f(x;\boldsymbol \theta)^q}{1-q}.
\end{equation}
\cite{Jackson08,Wadatwopara}.

If we consider the probability values from $f(x_i;\boldsymbol \theta)$ for random observations $x_1,x_2,\cdots,x_n$, then we have Tsallis $q$-entropy defined as 
\begin{equation}\label{Tsallisq}
-\sum_{i=1}^{n} \frac{f(x_i;\boldsymbol \theta)-f(x_i;\boldsymbol \theta)^q}{1-q}=-\frac{1-\sum_{i=1}^{n} f(x_i;\boldsymbol \theta)^q}{1-q}= -\sum_{i=1}^{n}f(x_i;\boldsymbol \theta)^q\log_q(f(x_i;\boldsymbol \theta)).
\end{equation}
$\log_q(f(x_i;\boldsymbol \theta)) = \frac{f(x_i;\boldsymbol \theta)^{1-q}-1}{1-q}$ is $q$-deformed form of ordinary logarithm. For $q \rightarrow 1$, $\log_{q}$ drops to $\log$ \cite{Tsallisbook09}.  $f(x;\boldsymbol{\theta})^q$ corresponds to (unnormalized) escort distribution connected with the $q$-calculus \cite{Ernst12}. Escort distribution, or ‘‘zooming distribution’’, was originally realized by  relation with dynamical chaotic systems \cite{Beck97,Beck04} and is also broadly applied into multifractals \cite{Harte01}. 

\section{Estimation methods based on maximum likelihood type}\label{Llocscale}
The maximum distorted likelihood estimation (MDLE)  is a tool to estimate the parameters $\boldsymbol \theta$ of model $f: x \times \boldsymbol \theta \in \mathbb{R}^p \rightarrow \mathbb{R}$ which is a probability density (p.d.) function corresponding to random variables $X_1$,\dots,$X_n$ which are non-identical and so some of $X_i$ are disorder or a member of another  parametric model $f(x;\boldsymbol{\tau})$. The likelihood function $L^D$ is given by
\begin{equation}\label{aplusformula1}
L^D(f(\boldsymbol{x};\boldsymbol{\theta})) = \prod_{i=1}^{n} [\beta + f(x_i;\boldsymbol \theta)],~~\beta \geq 0,
\end{equation}
$x_i$ in a vector $\boldsymbol{x}=(x_1,x_2,\cdots,x_n)$ is the observed value of random variable $X_i$. Eq. \eqref{aplusformula1} is maximized according to parameters $\boldsymbol \theta$ in order to get estimators $\hat{\boldsymbol{\theta}}$  when random variables $X_i$ are independent and non-identically distributed. When $\beta=0$, MDLE drops to MLE in which random variables are assumed to be identical. The tuning constant (TC) $\beta$  makes a perturbation to $f(x_i;\boldsymbol \theta)$.  This perturbation can be a good tricker to perform an efficient modeling when observations $x_1,x_2,\cdots,x_n$ in a data set are distributed non-identically \cite{Vajda86}. The maximization of likelihood function in Eq.  \eqref{aplusformula1} with $\beta=0$ coincides to minimization of entropy in Eq. \eqref{Tsallisq} with $q=1$ \cite{PardoSD,AbeOka01book,KullbackIT}, because minimization between $f(x;\boldsymbol \theta)$ and $f(x;\hat{\boldsymbol{\theta}})$ is performed. Further, Fisher information can be derived by use of entropy minimization \cite{CanKor18,Plastinoetal97}. 

\section{Generalized maximum likelihood estimation: M-estimations for parameters of exponential power distribution}\label{prelimi}

A generalized maximum likelihood estimation is defined by taking $\Lambda$ of likelihood $L$ in Eq. \eqref{aplusformula1}, i.e. $\Lambda$-likelihood, which is for $n$ case. If we consider a case in which $n=1$, then we will have a functional form of $\Lambda(f(x;\boldsymbol \theta))$. Thus, we have M-functions \cite{CanKor18} which are obtained by the generalized version of likelihood function with $n=1$ representing the objective function given by  \cite{Hub64}
\begin{equation}\label{defofrho}
\rho(x;\boldsymbol \theta)=\Lambda(f(x;\boldsymbol \theta)), 
\end{equation}
where $\Lambda:\mathbb{R} \rightarrow \mathbb{R}$ is a concave (or convex) function which can be taken as $\log$ and $\log_q$  \cite{CanKor18}. Then, objective functions based on $\log$ and $\log_q$ are obtained. A vector $\Psi$ for score functions $\psi_1,\psi_2,\cdots,\psi_p$ is given by
\begin{equation}\label{defofpsi}
\Psi(x; \boldsymbol \theta) = \frac{\partial}{\partial \boldsymbol \theta} \rho(x; \boldsymbol \theta),
\end{equation}
\noindent where vector $\boldsymbol \theta \in \mathbb{R}^p $ represents the number of $p$ parameters $\theta_1, \theta_2$ and $\theta_p$, $p \in \mathbb{Z}^{+}$ from a parametric model $f(x;\boldsymbol \theta) \in \mathbb{R}$.

M-estimations are defined by an objective function given as $\rho: x  ~ \times ~ \boldsymbol \theta \in \mathbb{R}^p \rightarrow \mathbb{R}$ for the number $n$ of observations represented by $x_1,x_2,\cdots,x_n$, minimizing 
\begin{equation}\label{defofrhoF}
\sum_{i=1}^n \rho(x_i; \boldsymbol \theta)_{|\boldsymbol \theta=\hat{\boldsymbol \theta}},
\end{equation}
over $\hat{\boldsymbol \theta}$ or by a vector of score functions given as $\Psi: x  ~ \times ~ \boldsymbol \theta \in \mathbb{R}^p \rightarrow \mathbb{R}^p$ as  solution for $\hat{\boldsymbol \theta}=(\hat{\theta}_1,\hat{\theta}_2,\cdots,\hat{\theta}_p)$ from a system of estimating equations 
\begin{equation}\label{defofpsiF}
\sum_{i=1}^n \Psi(x_i; \boldsymbol \theta)_{|\boldsymbol \theta=\hat{\boldsymbol \theta}}=\boldsymbol 0,
\end{equation}
which are used to get M-estimators $\hat{\boldsymbol \theta}$ from M-functions based on $\rho$ and $\psi$ in Eqs. \eqref{defofrho} and \eqref{defofpsi}, respectively. If $\Lambda$ is $\log$ and $\log_q$, then $\Psi$ is  $\log$ and $\log_q$ score function from derivatives of $\rho$ w.r.t $\boldsymbol{\theta}$, respectively \cite{Hub64,God60,Hampeletal86}. 

\subsection{Estimating equations based on $\log$ likelihood for parameters of exponential power distribution}\label{introEEs}
Gamma distribution is solution of ordinary differential equation \cite{SmithODE}. The continuous $Y \in [0,\infty)$ and discrete $Z=\{1,-1\}$ are variables distributed as gamma and probability $1/2$ around location $\mu$, respectively.  If we apply variable transformation $X=Y^{1/\alpha}Z$ to gamma distribution with special values of parameters, we  produce exponential power (EP) distribution which is a flexible function to model shape of data sets. The properties such as existence of moments and unibomodality are important indicators to use this function for modeling data sets (see \cite{CanEnt18}). The p.d. function of EP is given by  
\begin{equation}\label{etaESEPpdf}
f(x;\mu,\sigma,\alpha) = \frac{\alpha}{2 \sigma \Gamma(1/\alpha) } \exp \{- \left( \frac{|x - \mu|}{\sigma} \right)^\alpha   \}, ~~ x, \mu \in \mathbb{R}, \sigma >0,  \alpha>0.
\end{equation}
\noindent $\mu$ and $\sigma$ are location and scale parameters, respectively. The parameter $\alpha$ controls shape and peakedness of function $f$ and it is important to increase the modeling capability of $f$.  

MLE is used to produce the estimating equations (EEs). The parameters $\mu$, $\sigma$ and $\alpha$ in EP are estimated by use of EEs which generate score function $S$. In order to get EEs for these parameters, log-likelihood function of EP distribution in Eq. \eqref{etaESEPpdf} is necessary and it is given by 
\begin{equation}\label{LogLSpecialform2Can}
\log[L(\boldsymbol{x};\mu,\sigma,\alpha)]=n \log \bigg( \frac{\alpha}{2 \sigma \Gamma(1/\alpha)} \bigg) - \sum_{i=1}^{n} \left( \frac{|x_i - \mu|}{\sigma} \right)^\alpha.
\end{equation}
The derivatives of $\log(L)$ in Eq. (\ref{LogLSpecialform2Can}) with respect to (w.r.t) parameters $\mu$, $\sigma$ and $\alpha$ are taken and setting them to zero in order to obtain EEs. After algebraic manipulations are performed, we get the following EEs based on $\log$-score function $S$ for the parameters $\mu$, $\sigma$ and $\alpha$:
\begin{equation}\label{hatmuEE}
\hat{\mu}= \frac{\sum_{i=1}^{n} m(x_i;\hat{\mu},\hat{\sigma},\hat{\alpha}) x_i }{\sum_{i=1}^{n}  m(x_i;\hat{\mu},\hat{\sigma},\hat{\alpha})},
\end{equation}
\begin{equation}\label{hatsigmaEE}
\hat{\sigma} =\bigg[ \frac{1}{n}\sum_{i=1}^{n} m(x_i;\hat{\mu},\hat{\sigma},\hat{\alpha}) (x_i - \hat{\mu})^2 \bigg]^{1/2},
\end{equation}
\begin{equation}\label{hatalphaEE}
\hat{\alpha} = \left[ \sum_{i=1}^n \left[n \left(1+\frac{\psi(1/\hat{\alpha})}{\hat{\alpha}} \right)  \right]^{-1} \left\{\frac{m(x_i;\hat{\mu},\hat{\sigma},\hat{\alpha})}{\hat{\alpha}} \left( \frac{x_i-\hat{\mu}}{\hat{\sigma}} \right)^2 \log\left( \frac{|x_i- \hat{\mu}|}{ \hat{\sigma}}  \right)  \right\} \right]^{-1},
\end{equation}
The estimators $\hat{\mu}$, $\hat{\sigma}$ and $\hat{\alpha}$ are MLE of parameters $\mu$, $\sigma$ and $\alpha$. Here,
\begin{equation}\label{weight.func.genMLE}
m(x_i;\hat{\mu},\hat{\sigma},\hat{\alpha})=\hat{\alpha} \left(\frac{|x_i - \hat{\mu}|}{\hat{\sigma}}\right)^{\hat{\alpha}-2}=
\hat{\alpha}|y_i|^{\hat{\alpha}-2}
\end{equation}
\noindent  is a common function for MLE of parameters $\mu$, $\sigma$ and $\alpha$. It is also known as a weight function \cite{Cantez19}. $\psi$ in Eq. \eqref{hatalphaEE} is digamma function. The equations in (\ref{hatmuEE})-(\ref{hatalphaEE}) are EEs of parameters. For $\alpha=2$ known as Gaussian distribution, EEs of $\mu$ and $\sigma$ are given by \cite{MLqEBio}. The functions  depending only $x$, $\mu$ and $\sigma$ in Eqs. (\ref{hatmuEE})-(\ref{hatsigmaEE})  can be constructed  for the estimations of parameters $\mu$ and $\sigma$ if we have a location-scale model for an arbitrary function.  

\subsection{Structure of Huber's M-functions from EP distribution}

Let $S$ be a log-score function of EEs from $\log(f(x;\boldsymbol \theta))$. Then, $S$ can be obtained if it is written as   following form:
\begin{equation}\label{scorelogf}
S(y) = m(y) \cdot |y|,
\end{equation}
\noindent where $y=\frac{x - \mu}{\sigma}$. The common function $m$ for two parameters $\mu$ and $\sigma$ has also been defined by Eq. (\ref{weight.func.genMLE}) as a closed form of these parameters. $\rho(y)=-\log(f)$ is an objective function to show a connection between the getting log-score function $S$ produced by EEs and the derivative of $\rho(y)$ w.r.t the variable $y$. Thus the correspondence between maximum log-likelihood estimation method used to derive objective function $\rho$ and EEs which produce the score function for location and scale parameters are clarified.    The derivative of objective function $\rho$ w.r.t $y$ will be a log-score from function $m$, because if $\rho(y)=|y|^{\alpha}$, then $\frac{d}{dy}\rho(y)=S(y)=\alpha |y|^{\alpha-1}\text{sign}(y)$ from Eq. (\ref{scorelogf}) for EP distribution  (see p.52 in \cite{Maronna76} and p. 5 in \cite{Bercher12}). Thus we can have same mathematical expression with Eq.  (\ref{scorelogf}). The function $S$ defined as a log-score function is in fact a member of EEs and is used for estimation instead of using $\rho$ \cite{God60,PhDthesis,GodTh84,ProfMary}. 

The objective function of Huber is defined as 
\begin{equation}\label{Huberrho}
\rho(y)=
\left\{
\begin{array}{ll}
y^2~ &,  |y| \leq r; \\
2r|y|-r^2 &, |y| > r,
\end{array}
\right.
\end{equation}
\noindent and score function $\frac{d}{dy} \rho(y)=2S^H(y)$ is given by
\begin{equation}\label{Huberpsi}
S^H(y)=
\left\{
\begin{array}{ll}
y &, |y| \leq r; \\
\text{sign}(y)r &, |y| > r.
\end{array}
\right.
\end{equation}
Huber's score function in Eq. \eqref{Huberpsi} is a combination of log-score functions from normal and Laplace for $\alpha=2$ and $\alpha=1$ in EP distribution respectively from Eq. (\ref{etaESEPpdf}) \cite{Hub64}. The variable parts of Eq. \eqref{Huberpsi} are obtained by taking derivative w.r.t variable $y$ of $\rho(y)$, as given by Eq. \eqref{scorelogf} as base of  Eq. \eqref{Huberpsi} produced originally from spirit of EEs. Since  $S^H$ in Eq. \eqref{Huberpsi} is used with EEs,  $S^H$ with EEs is known as M-estimation.

\subsection{M-estimators $\hat{\mu}$ and $\hat{\sigma}$ in estimating equations based on combined $\log$-score functions  for non-identically distributed case of  EP distribution} \label{Riemannintegration}
The combination of log-score function $S$ with different values of parameter $\alpha$ in EP will be used to get M-estimators $\hat{\mu}$ and $\hat{\sigma}$  which are efficient for parameters $\mu$ and $\sigma$, because we want to propose two generalizations of Huber's score function from EEs. Inspiring from  Huber's M-functions \cite{Hub64} and using the spirit of EEs \cite{God60}, new score functions $S^{*}$ and $S^{*H}$ will be proposed by use of EEs in Eqs. \eqref{hatmuEE} and \eqref{hatsigmaEE} for the parameters $\mu$ and $\sigma$, respectively. After $S$ in Eq. \eqref{scorelogf} is partitioned, we have partial form of $S$. As a result, we have $\alpha_1$, $\alpha_2$ and $\alpha_3$ for left, middle and right parts on the real line, as proposed by Huber's $S^H$ in Eq. \eqref{Huberpsi} which is mainly from EEs in \cite{God60}. $n_1$, $n_2$ and $n_3$ represent sample size for these three parts. Each part $j$ has its corresponding score functions $S_j$.  The estimators $\hat{\mu}$ and $\hat{\sigma}$ from EEs are given by the following forms, respectively.

\begin{equation}\label{MLEmnmu}
\hat{\mu} = \sum_{j=1}^{3} \sum_{i=1}^{n_j} S_j(x_{ji};\hat{\mu},\hat{\sigma},\alpha_j)|y_{ji}|^{-1}x_{ji} / \sum_{j=1}^{3} \sum_{i=1}^{n_j} S_j(x_{ji};\hat{\mu},\hat{\sigma},\alpha_j)|y_{ji}|^{-1},
\end{equation}
\begin{equation}\label{MLEmnsigma}
\hat{\sigma} =  \left[ \frac{1}{\sum_{j=1}^{3} n_j} \sum_{j=1}^{3} \sum_{i=1}^{n_j} S_j(x_{ji};\hat{\mu},\hat{\sigma},\alpha_j)|y_{ji}|^{-1}(x_{ji} - \hat{\mu})^2 \right]^{1/2},
\end{equation}
where $y_{ji}=\frac{x_{ji} - \hat{\mu}}{\hat{\sigma}}$, $k,~t>0$ and we have
\begin{equation}\label{Scorea1a2a3a4nosign}
S^{*}(y)=
\left\{
\begin{array}{ll}
S_1 =  \alpha_1 |y|^{\alpha_1-1} &, (-\infty,-k); \\
S_2 =  \alpha_2 |y|^{\alpha_2-1} &, [-k,t]; \\
S_3 =  \alpha_3 |y|^{\alpha_3-1} &, (t,\infty).
\end{array}
\right.
\end{equation}
Huber's M-function in Eq. \eqref{Huberpsi} is generalized by partial log-score function $S^{*}$ in Eq. (\ref{Scorea1a2a3a4nosign}) which should not be continuous at the points $-k$ and $t$ for $\alpha_1, \alpha_2$ and $\alpha_3$, because modeling data, i.e. Riemann integration's rule is applied and histograms on the real line are constructed randomly, is equivalent to an integration of area under function. As it is well-known,  integration works even though the function $S^{*}$ used to fit data set is not continuous \cite{CanEnt18,LehmannCas98}. Thus, the discontinuity from $S_1$ to $S_3$ is not a problem for estimations of parameters. $S^{*}$ is log-score function to get M-estimators $\hat{\mu}$ and $\hat{\sigma}$ from EEs. $S^{*}$ can be replaced with $S^{*H}$, or alternatively, one can use  different score functions in order to manage the efficiency. Thus, we have an efficient fitting of data.

We propose a log-score function $S^{*H}$ given by
\begin{equation}\label{Scorea1a2a3a4Huber}
S^{*H}(y)=
\left\{
\begin{array}{ll}
S_1 = -k \alpha_1 |y|^{\alpha_1-1} &, (-\infty,-k); \\
S_2 =  \alpha_2 |y|^{\alpha_2-1} &, [-k,t]; \\
S_3 =  t  \alpha_3 |y|^{\alpha_3-1} &, (t,\infty),
\end{array}
\right.
\end{equation}
which is proposed by use of Huber's score function in Eq. \eqref{Huberpsi}. The log-score function $S^{*H}$ in Eq. \eqref{Scorea1a2a3a4Huber} is a generalized form of Eq. \eqref{Huberpsi}. Since the shape parameters $\alpha_1$, $\alpha_2$ and $\alpha_3$ in Eq. \eqref{Scorea1a2a3a4nosign} and (\ref{Scorea1a2a3a4Huber}) are added, we will have robust and efficient M-estimators for $\mu$ and $\sigma$ when we compared with Huber's $S^H$ which is not capable to fit shape of data sets. We also provide that Huber's $S^H$ in Eq. \eqref{Huberpsi} works on the base of EEs. Note that $r$ of Huber's M-functions from Eqs. \eqref{Huberrho} and \eqref{Huberpsi} can be replaced by $-k$ and $t$ because of spirit of EEs and estimations of paramaters \cite{LehmannCas98,CanEnt18} from Riemann integration.

\subsection{Estimating equations based on $\log_q$ and distorted $\log$ likelihoods for parameters $\mu$, $\sigma$ and $\alpha$  of EP distribution} \label{Lqlocscale}
As an alternative tool for robust estimation, the generalized and deformed entropies/logarithms \cite{GentEntReview1,JizbaArimitsu,JizKorhybrid,Wadatwopara,Jan17,Jan18,Jan19} and the divergences for estimation \cite{Basuetal98,Daroczy} can be used if the observations $x_1,x_2,\cdots,x_n$ in a data set are distributed non-identically. We will use $q$-deformed logarithm for the sake of its concavity, analytical simplicity and advantage for the constrained optimization in computation for parameter $q$ \cite{FerrariYang10,CanKor18,MLqEBio,FerrariOpere}. In MLE, it is possible to take the logarithm of likelihood function for a simple calculation in order to get estimators of parameters. In such a situation, log can be replaced with $\log_q$. Then it is known as maximum $\log_q$ likelihood estimation (MqLE) \cite{FerrariYang10,MLqEBio,Vajda86}. 

After taking $\log$ for Eq. \eqref{aplusformula1} with $\beta=0$ and replacing $\log$ form of Eq. \eqref{aplusformula1} with $\log_q$ and $\beta=0$, MqLE for parameters of EP distribution are obtained by optimizing the following function according to parameters  $\mu, \sigma$ and $\alpha$
\begin{equation}\label{qMLEmn}
l_q (L(\boldsymbol{x};\mu, \sigma, \alpha)) = \sum_{i=1}^{n} \log_q[f(x_i;\mu, \sigma, \alpha)],
\end{equation}
$\boldsymbol{x}$ is a vector of random observations $x_1,x_2,\cdots,x_n$. The p.d. function $f$ is in Eq. \eqref{etaESEPpdf}. The derivatives of $l_q(L)$ in Eq. (\ref{qMLEmn}) w.r.t parameters $\mu$, $\sigma$ and $\alpha$ are taken and setting them to zero in order to obtain EEs. After algebraic manipulations are performed, we get the following EEs based on $\log_q$-score function $S_q$ for parameters $\mu$, $\sigma$ and $\alpha$:
\begin{equation}\label{qMLEmnmu}
\hat{\mu}=\sum_{i=1}^{n}S_q(x_i;\hat{\mu},\hat{\sigma},\hat{\alpha})|y_i|^{-1}x_i/\sum_{i=1}^{n} S_q(x_i;\hat{\mu},\hat{\sigma},\hat{\alpha})|y_i|^{-1},
\end{equation}
\begin{equation}\label{qMLEmnsigma}
\hat{\sigma} =  \left[\frac{1}{\sum_{i=1}^{n}w_i}\sum_{i=1}^{n}S_q(x_i;\hat{\mu},\hat{\sigma},\hat{\alpha})|y_i|^{-1}(x_i-\hat{\mu})^2 \right]^{1/2},
\end{equation}
\begin{equation}\label{qMLEmnalpha}
\hat{\alpha}=\hat{\alpha}\frac{\sum_{i=1}^{n} S_q(x_i;\hat{\mu},\hat{\sigma},\hat{\alpha})|y_i|\log(|y_i|)}{\sum_{i=1}^{n} w_i}-\psi(1/\hat{\alpha}),
\end{equation}
where $S_q=w_i S$ is a weighted ($w$) form of log-score function $S$, $w_i=f(x_i;\hat{\mu},\hat{\sigma},\hat{\alpha})^{1-q},~ S(x_i;\hat{\mu},\hat{\sigma},\hat{\alpha})=\hat{\alpha} |y_i|^{\hat{\alpha}-1}$, $y_i = \frac{x_i-\hat{\mu}}{\hat{\sigma}}$. $\psi$ is digamma function. $S_q$ from EEs can be used to fit a non-identically distributed data set.  For $q=1$, $S_q$ drops to $S$ in Eq. \eqref{scorelogf}. $S$ can be derived by a function without normalizing factor, such as escort distribution \cite{Jizba3} or $S$ can be an arbitrary function in M-functions  \cite{God60}. 

For MDLE of the parameters of $f(x;\boldsymbol \theta)$, logarithmic form of Eq. \eqref{aplusformula1}, i.e. $\log(\beta+f)$, is used. $S_q$ in Eqs. \eqref{qMLEmnmu}-\eqref{qMLEmnalpha} can be replaced with the distorted log-score function $S^D=w_iS$ and $w_i$ is $\frac{f(x_i;\hat{\mu},\hat{\sigma},\hat{\alpha})}{\beta+f(x_i;\hat{\mu},\hat{\sigma},\hat{\alpha})}$; because, one get  same result, which shows us an interesting connection between  $\log_q$ and distorted log-likelihoods for parameters of EP distribution. Note that $S_q$ and $S^D$ are comparable due to the weighted form of $\log$-score function $S$, i.e. $wS$ which can produce robust and also efficient M-estimators $\hat{\boldsymbol{\theta}}$. $S_q$ and  $S^D$ are redescending M-functions which are used with EEs to get M-estimators \cite{Hampeletal86}, because $S_q$ and $S^D$ are zero if $\lim_{x \rightarrow \infty} f(x;\mu,\sigma,\alpha)=0$.  It can be necessary to have a partial form of $S_q$ or $S^D$ which is used to model a data set more efficiently for the estimations of $\mu$ and $\sigma$ when many outliers or huge disorders exist in a data set. 

\subsection{Robustness of M-estimators: Score functions of parameters $\mu$, $\sigma$ and $\alpha$ of EP distribution for $\log_q$, $\log$ and distorted $\log$ functions}\label{robustscore}
If we observe the non-identically distributed random variables $X_1,X_2,\cdots,X_n$, we consult the robust statistics which can represent the behavior of bulk of data in a data set when there exist outlier(s) in the data set.  In order to test the robustness of M-estimators produced by objective $\rho$ and score $\Psi$ functions used to model a data set with outlier(s),  we need to examine the values of score functions when $y$ goes to infinity.  \cite{Hub64,Hampeletal86}. 

Let us examine whether or not the elements of a vector $\Psi$ for parameters of EP distribution in $\log$, $\log_q$ and distorted log functions are finite. For the sake of simplicity, the finiteness for positive side of score functions in the vector $\Psi$ is examined, because the distribution is symmetric \cite{CanEnt18}. If all elements of vector $\Psi$ which shows robustness property of M-estimators are finite for $y \rightarrow \infty$, then M-estimators will be robust \cite{Hampeletal86}.

Let $\Psi_{\log_q}=(^\mu \psi_{\log_q},^\sigma \psi_{\log_q},^\alpha \psi_{\log_q})$ be a vector of score functions from derivative of $\Lambda(f(x;\boldsymbol{\theta}))$ w.r.t parameters $\boldsymbol{\theta}=(\mu,\sigma,\alpha)$ of EP distribution. Here $\Lambda$ can be $\log$ and $\log_q$. Let us examine the finiteness of score functions.
\begin{equation}\label{logqM}
^\mu \psi_{\log_q}(x;\boldsymbol{\theta})=\frac{\partial}{\partial \mu} \log_q \left(f(x;\boldsymbol{\theta}) \right) =-f(x;\boldsymbol{\theta})^{1-q} \alpha  y^{\alpha-1},
\end{equation}
\begin{equation}\label{logqS}
^\sigma \psi_{\log_q}(x;\boldsymbol{\theta})=\frac{\partial}{\partial \sigma} \log_q \left(f(x;\boldsymbol{\theta}) \right)  = f(x;\boldsymbol{\theta})^{1-q} \left[ -\frac{1}{\sigma} + \frac{\alpha}{\sigma} y^\alpha \right],
\end{equation}
\begin{equation}\label{logqA}
^\alpha \psi_{\log_q}(x;\boldsymbol{\theta})=\frac{\partial}{\partial \alpha} \log_q \left(f(x;\boldsymbol{\theta}) \right)  = f(x;\boldsymbol{\theta})^{1-q} \left[ \frac{1}{\alpha} + \frac{\psi(1/\alpha)}{\alpha^2} - y^\alpha \log(y) \right],
\end{equation}
\input{table1.tex}
where $y=\frac{x-\mu}{\sigma}$. 

Let $\Psi_{\log}^D=(^\mu \psi_{\log}^D,^\sigma \psi_{\log}^D,^\alpha \psi_{\log}^D)$ be a vector of score functions from derivative of $\log(\beta+f(x;\boldsymbol{\theta}))$ w.r.t parameters of EP distribution. Let us examine the finiteness of score functions.
\begin{equation}\label{logDM}
^\mu \psi_{\log}^D(x;\boldsymbol{\theta})=\frac{\partial}{\partial \mu} \log \left(\frac{f(x;\boldsymbol{\theta})}{\beta+f(x;\boldsymbol{\theta})} \right) =- \frac{f(x;\boldsymbol{\theta})}{\beta+f(x;\boldsymbol{\theta})} \alpha  y^{\alpha-1},
\end{equation}
\begin{equation}\label{logDS}
^\sigma \psi_{\log}^D(x;\boldsymbol{\theta})=\frac{\partial}{\partial \sigma} \log \left(\frac{f(x;\boldsymbol{\theta})}{\beta+f(x;\boldsymbol{\theta})} \right) = \frac{f(x;\boldsymbol{\theta})}{\beta+f(x;\boldsymbol{\theta})}  \left[ -\frac{1}{\sigma} + \frac{\alpha}{\sigma} y^\alpha \right],
\end{equation}
\begin{equation}\label{logDA}
^\alpha \psi_{\log}^D(x;\boldsymbol{\theta})=\frac{\partial}{\partial \alpha} \log  \left(\frac{f(x;\boldsymbol{\theta})}{\beta+f(x;\boldsymbol{\theta})} \right) =\frac{f(x;\boldsymbol{\theta})}{\beta+f(x;\boldsymbol{\theta})}  \left[ \frac{1}{\alpha} + \frac{\psi(1/\alpha)}{\alpha^2} - y^\alpha \log(y) \right],
\end{equation}
where $y=\frac{x-\mu}{\sigma}$.
\input{table2.tex}

M-estimators $\hat{\mu}$, $\hat{\sigma}$ and $\hat{\alpha}$ will be robust if $0 < q <1$ and $\beta>0$, because we have finite values of limit  for score functions (see Tables \ref{qvaluetable} and \ref{Dvaluetable}). Since $y$ is a variable which includes location and scale parameters with variable $x$, i.e.  $y=\frac{x-\mu}{\sigma}$, examining the finiteness of variables $y$ and $x$ which go to infinity in $\Psi$ is equivalent to each other \cite{Hampeletal86,LehmannCas98}.

\section{Fisher information matrices based on derivative of score functions  and information theory}\label{FisherEEs}
Fisher information (FI) is used to get variance-covariance of MLE. FI matrix of objective function $\rho=-\log(f(x;\boldsymbol{\theta}))$, that is $-E\left[\frac{\partial^2}{\partial \boldsymbol{\theta} \partial \boldsymbol{\theta}^T}\log[f(x;\boldsymbol{\theta})]\right]=E\left[\frac{\partial^2}{\partial \boldsymbol{\theta} \partial \boldsymbol{\theta}^T}\rho[f(x;\boldsymbol{\theta})]\right]=E\left[\left\{\frac{\partial}{\partial \boldsymbol{\theta} }\rho[f(x;\boldsymbol{\theta})]\right\}\left\{\frac{\partial}{\partial \boldsymbol{\theta} }\rho[f(x;\boldsymbol{\theta})]\right\}^T\right]$ as the known definition of FI, \cite{AmariIG,Fisher25} can be adopted for score functions $S^H$, $S^*$, $S^{*H}$, $S_q$ and $S^D$ from EEs  by use of the definition of Tsallis $q$-entropy as a generalization of Shannon entropy \cite{CanKor18,Plastinoetal97,CanKor16}. Thus, FI of score functions based on $S$ will be proposed for M-estimators as a generalization of MLE (see \ref{proofofFisher}). Since FI matrices are based on derivative of score functions w.r.t parameters $\boldsymbol \theta$, they are comparable \cite{CanKor18,infocomplestamod,MyungBalasubramanian,ClusteringFisher,Maybankvol}. Note that since M-functions can become $S(x;\boldsymbol \theta)$ (or $S$ can be an arbitrary function) which can be generated from $\Psi(x;\boldsymbol \theta)$, it is necessary to use FI based on  the derivative of score functions w.r.t $\boldsymbol{\theta}$ instead of using variance-covariance of M-estimators produced by M-functions from $\rho$ and $\psi$ \cite{Hub81,GodHeyquasi}.

FI matrix obtained from derivative of log-score function $S^*$ in Eq. \eqref{Scorea1a2a3a4nosign} or $S^{*H}$ in Eq. \eqref{Scorea1a2a3a4Huber}  w.r.t parameters $\mu$ and $\sigma$ is defined as 
\begin{equation}\label{derivSFisher}
F_{\log}(S_j;f_2;\boldsymbol \theta,\alpha_j)=  n \sum_{j=1}^3 \int_{a_j}^{b_j} \frac{\partial}{\partial \boldsymbol \theta} S_j(x_j;\boldsymbol \theta,\alpha_j) \frac{\partial}{\partial \boldsymbol \theta} S_j(x_j;\boldsymbol \theta,\alpha_j)^T f_2(x_j;\boldsymbol \theta,\alpha_2) \ud x_j,
\end{equation}
where $\boldsymbol \theta=(\mu,\sigma)$,~ $a_j,b_j \in \mathbb{R}$ and $S_j=\alpha_j(\frac{|x_j-\mu|}{\sigma})^{\alpha_j-1}\text{sign}(x_j-\mu)$.

It is possible to propose FI for location, scale and shape parameters if we use distributions which can derive score functions via using MLE, MqLE and MDLE for these parameters. FI matrix obtained from $\log_q$-score function $S_q$ is defined as
\begin{equation}\label{derivSFisherq}
F_{\log_q}(S_q;f;\boldsymbol \theta)= n  \int_{a}^{b} [(1-q)S^2(x;\boldsymbol \theta) \frac{\partial}{\partial \boldsymbol \theta} S(x;\boldsymbol \theta)^T+\frac{\partial}{\partial \boldsymbol \theta}S(x;\boldsymbol \theta)\frac{\partial}{\partial \boldsymbol \theta}S(x;\boldsymbol \theta)^T] f(x;\boldsymbol \theta)^{2-q} \ud x.
\end{equation}

FI matrix obtained from distorted log-score function $S^D$ is defined as
\begin{footnotesize}
\begin{equation}\label{derivSFisherbeta}
F_{\log}(S^D;f;\boldsymbol \theta)= n  \int_{a}^{b} \left[\left(\frac{\beta}{\beta+f(x;\boldsymbol \theta)}\right)S^2(x;\boldsymbol \theta) \frac{\partial}{\partial \boldsymbol \theta} S(x;\boldsymbol \theta)^T  + \frac{\partial}{\partial \boldsymbol \theta}S(x;\boldsymbol \theta)  \frac{\partial}{\partial \boldsymbol \theta}S(x;\boldsymbol \theta)^T \right]\frac{f(x;\boldsymbol \theta)^2}{\beta+f(x;\boldsymbol \theta)}  \ud x,
\end{equation}
\end{footnotesize}
where $\boldsymbol \theta=(\mu,\sigma,\alpha)$, ~ $a,~b \in \mathbb{R}$ and $S=\alpha(\frac{|x-\mu|}{\sigma})^{\alpha-1}\text{sign}(x-\mu)$. The inverse of $F_{\log_q}$ and its special case $F_{\log}$ are defined as variance-covariance matrices of M-estimators $\hat{\boldsymbol \theta}$. The variances-covariances of $\hat{\boldsymbol \theta}$ are denoted by $\text{Var}_{S_q}(\hat{\boldsymbol \theta})=F^{-1}_{\log_q}$ and $\text{Var}_{S}(\hat{\boldsymbol \theta})=F^{-1}_{\log}$ for $q=1$.  Eq. \eqref{derivSFisherq} with $q=1$ and Eq. \eqref{derivSFisherbeta} with $\beta=0$ are Eq. \eqref{derivSFisher} with parameter $\alpha=\alpha_j$.  The matrix in Eq. \eqref{derivSFisher} is positive semidefinite, i.e. non-negative, and symmetric. The matrices in Eqs. \eqref{derivSFisherq}-\eqref{derivSFisherbeta} are positive semidefinite and asymmetric if conditions are satisfied  (see  \ref{proofofFisherM} for details).  

\section{Selection of M-functions via tools from information theory}\label{volumeandIC}
\subsection{Selection of optimal M-function via volume} \label{modelselect}
The fitting performance of M-functions is tested by the volume of ellipsoid of M-function with underlying p.d. $f(x;\boldsymbol \theta)$. The volume based on score function $S$ can be proposed by replacing the known FI matrix based on objective function $\rho$ with FI matrix based on $S$ used in EEs. The main part of Eq. \eqref{entropygeneratorvol} is det$(F)$. Other parts are constants which do not affect volume exactly when volume based on $\rho$ in \cite{CanKor18} is compared with that in Eq.  \eqref{entropygeneratorvol}. Note that score functions derived by EEs are comparable for volume.

The volume is used to determine the values of $r,k,t$ in $S^H$, $S^{*}$ and $S^{*H}$ with EEs, $q$ and $\beta$ from $S_q$ and $S^D$ in MqLE and MDLE, respectively.  The different values of tuning constants (TCs) with estimated values of $\hat{\boldsymbol \theta}$ and fixed values of $\alpha$ are tried with grid search until the smallest values of volume and MAE are  obtained. Thus, the  optimal M-function is chosen among different values of TCs and the corresponding estimates according to TCs in M-functions. The volume is defined by 
\begin{equation}\label{entropygeneratorvol}
\text{Vol}_q(S_q;\hat{\boldsymbol \theta})=\left( \frac{2 \pi v}{n} \right)^{d/2} \frac{1}{\Gamma(d/2+1)} \frac{1}{\sqrt{\text{det}[F_{\log_q}(S_q;f;\hat{\boldsymbol \theta})]}},
\end{equation}
\noindent where $n$ is the number of observations in a data set. $v$ and $d$ are the number of eigenvalues and  dimension of FI matrix, respectively. $S_{q}$ is replaced with $S^H$, $S^{*}$, $S^{*H}$ and $S^D$ in $\log$ case to produce volumes based on $S^H$, $S^{*}$, $S^{*H}$ and  $S^D$, respectively. We have volume known as a general geometric measure for $d \geq 1$. In special case, if $\hat{\boldsymbol{\theta}}$ is a vector for $(\hat{\mu},\hat{\sigma})$ and  $(\hat{\mu},\hat{\sigma},\hat{\alpha})$, then dimensions of FI matrix are $d=2$ and $d=3$, respectively  \cite{CanKor18,infocomplestamod,MyungBalasubramanian,ClusteringFisher,Maybankvol,AmariIG}. 

\subsection{Selection of the best M-function via information criteria}\label{modelselectionICs} 
Information criteria (IC) are used to test whether or not the best M-function for a data set can be chosen. Since M-estimators are based on score functions  derived from $\log(f)$, $\log_q(f)$ and $\log(\beta+f)$, it is appropriate to use Akaike IC (AIC), corrected AIC (cAIC) and Bayesian IC (BIC) which should be based on score functions in EEs. Note that the known definitions of IC are comparable among different p.d. functions used to fit a data set, because p.d. functions are on the interval $[0,1]$. The ordinary logarithm from definition of IC is  used. So $\log(f(x;\boldsymbol \theta))$ will be comparable when IC are used for the selection of the best p.d. $f(x;\boldsymbol \theta)$ among p.d. functions.  In our case, IC are proposed by use of score functions $S$ instead of using objective functions $\log(L)$  defined as origin of IC. Since summation of a score function evaluated by values of $x_1,x_2,\cdots,x_n$ and $\hat{\boldsymbol{\theta}}$ is approximately zero and theoretically equals to zero, i.e. integration over all values of $x$ and true values of parameters \cite{AmariIG}, the absolute value of score function is taken in order to test the fitting performance of score function in its self \cite{CanKor18,AmariIG}. 

IC of the combined log-score functions is defined as 
\begin{equation}\label{AICformula}
IC_S=2 \sum_{j=1}^3  \sum_{i=1}^{n_j} |S_j(x_{ji};\hat{\mu},\hat{\sigma},\alpha_j)|  + c(p,n),
\end{equation}
where $c(p,n)$ is the given constant. Similarly, IC for $S_q$ is obtained by replacing $\log$ in Eq. (\ref{AICformula}) with $\log_q$. However, $S_q$ is not a partial form and so IC of $\log_q$-score function is defined as \cite{CanKor18,FerrariOpere}
\begin{equation}\label{AICformulalogq}
IC_{S_q}=2 \sum_{i=1}^n |S_q(x_i;\hat{\mu},\hat{\sigma},\hat{\alpha})| + c(p,n)
\end{equation}
for a determined value of $q$ via volume as well. $p$ is a number of the estimable parameters and $n$ is sample size. $S_q$ is replaced with $S^D$ for the distorted log likelihood.  The formulae in Eqs. \eqref{AICformula} and \eqref{AICformulalogq} are tools for selection of M-function. $r$, $-k$, $t$, $q$, $\beta$ and $\alpha$ in score functions are chosen for their corresponding score functions with appropriate minimum values of IC in Eqs. \eqref{AICformula} - \eqref{AICformulalogq} and the smallest value of MAE.  The performance of IC depends on the penalty term $c(p,n)$. When $c(p,n)$ is $2p$, $\frac{2pn}{n-p-1}$ and $p\log(n)$, penalty terms of AIC, cAIC and BIC are obtained, respectively \cite{Aka73}. Since penalty term of AIC is $2p$, the quality of AIC is low to select the best M-function. For this reason, cAIC and BIC are given. Their penalty terms are affected by $p$ and $n$ \cite{Hamparsum87,Roncpt97}. 

\section{Application for real data sets}\label{simrealapp}
Temperature measurements are important to observe the ecological movement and the results of temparature on earth such as sudden rains and floods should be examined precisely.  Grytviken temperature from  1905 to 2019 years is analyzed to get statistics about location, scale and also shape parameters of empirical distribution (Example 1 and see Table \ref{ex1temperature}) \cite{temperaturedata}.

Cancer treatment in genomic and pharmacology  includes measurements to examine affects of pharmacology  in treatments of genes.  Gene-drug correlation data should be analyzed precisely  to manage level of pharmacology (Example 2 and see Table \ref{ex2microarray})  \cite{microarray}. The real data set which is cDNA microarray coded as "SID W 486613, ESTs, Highly similar to Ovarian Granulosa Cell 13.0 KD protein HGR74 [Homo sapiens] [5':AA044350, 3':AA044028]"  
\cite{microarray} is used to get statistics as well.

If a distortion (occuring due to measurement error in experiment, difference in measurements from laboratory expert, incorrectly recorded data or the nature of  phenomena) in a data set exists, then the underlying distribution which represents the majority of a data set should be modelled efficiently. The location, scale and also shape parameters are mainly used to get the statistics about data analysis \cite{LehmannCas98,Hampeletal86}. Random number generation (see \ref{simdesigns}) is also an advantage to observe behaviours of real data sets. For these aims, we apply our EEs into estimation of these parameters robustly and precisely. The parameters $\alpha_1$, $\alpha_2$ and $\alpha_3$ are used to manage efficiency of score functions $S^{*}$ and $S^{*H}$ as well. The objective functions $\log_q(f)$ and $\log(\beta+f)$ are also applied to model real data sets. Since the modeling capability of the proposed functions for these data set is observed to be high, these real data sets have been chosen. Many phenomena can be modelled by these functions in order to be able to get efficient M-estimators from M-functions such as score and objective.   Two outlier values which are two times of maximal value of data set are added  as positive and negative values. Thus we added two outliers  into real data sets. The main aim is to observe the sensitivity of M-estimators from score functions with EEs and objective functions produced by MqLE and MDLE to two outliers.

\input{table3.tex}
\input{table4.tex}

Even if $S^{*}$ in Eq. \eqref{Scorea1a2a3a4nosign} or $S^{*H}$ in Eq. \eqref{Scorea1a2a3a4Huber} from EEs is infinite for $\alpha_1,\alpha_3>1$, the best score functions with EEs are used to model real data sets, which shows that the underlying and contamination should be efficiently modelled for observations from a finite sample size (see italic numbers in Tables \ref{sim1Design}-\ref{sim4Design}).  

Redescending M-functions derived by MqLE and MDLE can be alternative to each others for estimation and so we can have similar values for the estimation of parameters and MAE (see Table \ref{ex2microarray} for $S^D$, $S_q$, MDLE and MqLE). The value of MAE for $S^*$ is better than that of $S^{*H}$ and $S^H$ for case in Table 
\ref{ex2microarray}.  For overall assessment according to MAE, $S^D$  has  the smallest value of MAE among others in Table \ref{ex2microarray}. When we consider the results of MDLE in simulation in \ref{secsimulation} and real data sets (see for procedure in \ref{secreal}), MDLE should be preferred. 

\section{Conclusions and discussions}\label{labconclusions}
This paper has focused on derivation of score functions from EEs based on $\log$, $\log_q$ and distorted log-likelihoods and inference for location, scale and shape parameters of EP distribution. We have also proposed to estimate location and scale parameters for  shape parameter value determined by consulting FI from information geometry; IC which can be used as tools for model selection among competing models; and MAE used to make a comparison among M-functions all together.  We have considered to propose score functions as an approach in modeling. Even if we use EP distribution to produce score functions, one can use this approach for two aims. The first one is that score functions can be produced by use of different parametric models. The second one is that the score function in EEs is replaced by arbitrary ones which can model data set efficiently.

Applications from simulation and real data sets have been given to test the fitting performance of score functions with EEs and objective functions.   $\widehat{\text{MSE}}(\hat{\boldsymbol \theta})$ is used for the simulation in order to determine the values of TCs in score and objective functions. The results of simulation show that MDLE of parameters $\sigma$ and $\alpha$ of EP distribution are better than MqLE of these parameters. $q$ in MqLE affects the success of estimations of parameters $\sigma$ and $\alpha$.  We have showed that FI matrices in Eqs. \eqref{derivSFisher}, \eqref{derivSFisherq} and \eqref{derivSFisherbeta} are positive semidefinite if conditions are satisfied. Thus, the diagonal elements of inverse of FI matrix are variance of estimators, i.e. Var$(\hat{\mu})$, Var$(\hat{\sigma})$ and Var$(\hat{\alpha})$. Explicit expressions calculated for each element of FI for parameters $\mu$, $\sigma$ and $\alpha$ of EP distribution are given for $S^{*}$, $S^{*H}$ and $S_q$. 

M-estimators are robust for $q \in (0,1)$ and $\beta>0$. Robustness is not enough to imply the modeling capability for empirical distributions occurred by the underlying and contamination. The different score functions  with EEs should be used to model a data set even if the score functions $S^*$ and $S^{*H}$ derived by EEs are  infinite for parameters $\mu$ and $\sigma$ respectively when $\alpha>1$, as supported by the results of simulation for simultaneous estimation of parameters $\mu$ and $\sigma$. In addition, different objective functions should also be used to have the best fitting on a data set. In the application of real data sets, the role of robustness has been observed especially for estimation of $\sigma$ when outlier and non-outlier cases are compared. Redescending M-functions from $S_q$ and $S^D$ give efficient results extensively and they are finite for $q \in (0,1)$ and $\beta > 0$ when they are compared with score functions $S^{*H}$ and $S^*$ for estimations of $\mu$ and $\sigma$, as observed by simulation results. For determining the values of TCs, volume, IC and MAE should be used together.  

EEs for different p.d. functions will be derived.  EEs from estimation methods which are different from likelihood method, fractional objective functions, regression case and their multivariate forms with applications will be studied. FI matrix based on different estimation methods will also be derived to provide variance of the estimators obtained by their corresponding estimation methods. A user-friendly package  in open access statistical software $R$ will be prepared for practitioners in the applied science to manage precise modeling on the data sets.

\appendix

\section{Random number generation, tools for computation and simulation }\label{generationopti}
\subsection{Algorithm of random number generation}\label{simdesigns} 
The algorithm for number generation is given in the following order \cite{Elsal09}:
\begin{enumerate}
	\item Generate a random number $y$ from gamma distribution with shape parameter $1/\alpha$ and scale parameter $1$, i.e. $Y \sim \Gamma(1/\alpha,1)$,
	\item Use the variable transformation $x_j=\mu_j+ \sigma_j  Y^{1/\alpha_j}$, $j=1,2,3$, $x_2$ represents the artificial data set from the underlying distribution. $x_1$ and $x_3$ represent the artificial data sets from the contamination. Thus, the artificial random numbers are $\boldsymbol{x}=(x_1,x_2,x_3)$.
\end{enumerate}

\subsection{Tools used for computation}
All of computations are performed by using MATLAB R2013a.  The codes for genetic algorithm (GA) which is a derivative-free method to perform an optimization are given by the following order: 
\begin{enumerate}
	\item Provide lower and upper bounds for search space of GA: \\
	$\mu_{\text{lower}}=-10^{10}$;	$\mu_{\text{upper}}=10^{10}$;$\sigma_{\text{lower}}=0.01$;$\sigma_{\text{upper}}=10^{10}$;$\alpha_{\text{lower}}=0.01$;$\alpha_{\text{upper}}=10^{10}$;
	lb=[$\mu_{\text{lower}}$ $\sigma_{\text{lower}}$ $\alpha_{\text{lower}}$];ub=[$\mu_{\text{upper}}$ $\sigma_{\text{upper}}$ $\alpha_{\text{upper}}$];
	\item Use 'ga' to get the estimates of $\hat{\mu},\hat{\sigma},\hat{\alpha}$ by use of MqLE:  
	\begin{enumerate}
		\item 	opts = gaoptimset('CrossoverFcn',{@crossoversinglepoint},'display','off');
		\item 
		\begin{verbatim}1. function [Estimates] = MqLE(p,x,q)
		2. mu=p(1);sigma=p(2);alpha=p(3);
		3. f=(alpha/(2*sigma*gamma(1/alpha))).*exp(-(abs(x-mu)./sigma).^alpha);
		4. Lq=sum((f.^(1-q)-1)./(1-q)); 
		5. Estimates=-Lq;
		\end{verbatim}
		\item Give the estimates of  $\hat{\mu},\hat{\sigma},\hat{\alpha}$: \\	pMqLE=ga(@(p)MqLE(p,x,q),3,[],[],[],[],lb,ub,[],opts);
	\end{enumerate}
\end{enumerate}
The parameter values of GA are default of GA module in MATLAB R2013a. Since the optimization of objective functions in Eq. \eqref{aplusformula1} with $\log$, i.e. $\log(\beta+f)$, and Eq. \eqref{qMLEmn}  according to parameters is equivalent to the numerical solving of equations for the estimators $\hat{\mu}$, $\hat{\sigma}$ and $\hat{\alpha}$ in Eqs. \eqref{qMLEmnmu}-\eqref{qMLEmnalpha}, the objective functions $\log_q(f)$ and $\log(\beta+f)$ are used to get the estimates of $\hat{\mu}$, $\hat{\sigma}$ and $\hat{\alpha}$.  Note that $\alpha$ is an important parameter \cite{Orkcu15}. For this reason, GA is used. MDLE is obtained by same codes and  4. line is replaced with sum$(\log(\beta+f))$.

Iterative reweighting algorithm for simultaneous estimations of $\mu$ and $\sigma$ are given by \cite{Hub81}.  For computation of all EEs in Eqs. (\ref{MLEmnmu})-(\ref{MLEmnsigma}) and (\ref{qMLEmnmu})-(\ref{qMLEmnsigma}), initial values are provided by $\hat{\mu}=\text{Median}(\textbf{x})$, $\hat{\sigma}=\text{Median}(|x_i-\text{Median}(\textbf{x})|)$, $\boldsymbol{x}=\{x_1,x_2,\cdots,x_n\}$.

\subsection{The results of simulation for simultaneous estimation of $\mu, \sigma$ and simultaneous estimation of $\mu, \sigma, \alpha$}\label{secsimulation}
Simulation provides to observe the performance of different functions from score and objective. Thus, we make a comparison among their modeling capability for same types of designs in \ref{generationopti}. It is important  to drive a comprehensive simulation and its outputs are given by Tables \ref{sim1Design}-\ref{sim4Design} which show the results of simulation for M-estimates of $\hat{\mu}$, $\hat{\sigma}$ from EEs with score functions $S^H$, $S^{*H}$, $S^*$,  $S_q$, $S^D$ and M-estimates of $\hat{\mu}$, $\hat{\sigma}$, $\hat{\alpha}$ from MqLE and MDLE, the simulated variance ($\widehat{\text{Var}}(\hat{\boldsymbol \theta})$),  the simulated mean squared error ($\widehat{\text{MSE}}(\hat{\boldsymbol \theta})$) and TC values $k,t,q,\beta$ for $\hat{\mu}$, $\hat{\sigma}$ and $q,\beta$ for $\hat{\mu}, \hat{\sigma}$, $\hat{\alpha}$. TC$_1$, TC$_2$, TC$_3$ show values of TC for sample sizes $n=110,210,410$, respectively.  The values of TCs for robust and efficient estimations are determined until the smallest values of mean squared error (MSE) defined as $\frac{1}{m} \sum_{j=1}^{m} (\hat{\boldsymbol \theta}_j - \boldsymbol{\theta})^2$  are obtained.  If values of  $\widehat{\text{MSE}}(\hat{\boldsymbol \theta})$ and $\widehat{\text{Var}}(\hat{\boldsymbol \theta})$ are equal to each other, then biases of estimators $\hat{\boldsymbol \theta}$ are zero. 

We have four types of design. These designs will have outliers for $\alpha < 2$, because EP distribution becomes a heavy-tailed function if $\alpha \in (0,2)$. EP is more heavy-tailed function for $\alpha \leq 1$ \cite{Lucas97}. Thus, data which are distant from bulk of data can be generated. The sample sizes $n_1$ and $n_3$ of contamination are $5$, respectively. The sample size $n_2$ of underlying is $100,200,400$. The number of replication for sample sizes is $m=10^4$. The four designs  are given by following forms for Tables \ref{sim1Design}-\ref{sim4Design}, respectively:
\begin{itemize}
	\item $1^{th}$ design: $x_1=D_1(\alpha_1=1.1,\mu_1=5,\sigma_1=6,n_1)$, $x_2=D_2(\alpha_2=2,\mu_2=0,\sigma_2=1,n_2)$, $x_3=D_3(\alpha_3=1.2,\mu_3=4,\sigma_3=2,n_3)$.
	\item $2^{nd}$ design: $x_1=D_1(\alpha_1=1.1,\mu_1=2,\sigma_1=3,n_1)$, $x_2=D_2(\alpha_2=3,\mu_2=0,\sigma_2=1,n_2)$, $x_3=D_3(\alpha_3=1.2,\mu_3=3,\sigma_3=5,n_3)$.
	\item $3^{rd}$ design: $x_1=D_1(\alpha_1=1.2,\mu_1=3,\sigma_1=4,n_1)$, $x_2=D_2(\alpha_2=3,\mu_2=0,\sigma_2=1,n_2)$, $x_3=D_3(\alpha_3=0.8,\mu_3=3,\sigma_3=4,n_3)$.
	\item $4^{th}$ design: $x_1=D_1(\alpha_1=0.7,\mu_1=4,\sigma_1=2,n_1)$, $x_2=D_2(\alpha_2=1.3,\mu_2=0,\sigma_2=1,n_2)$, $x_3=D_3(\alpha_3=0.9,\mu_3=2,\sigma_3=3,n_3)$.
\end{itemize}
\input{table5.tex}
\input{table6.tex}
\input{table7.tex}
\input{table8.tex}
If  $\alpha_1,\alpha_3>1$, score functions $S^{*H}$ and $S^{*}$ are infinite when $y=\frac{x-\hat{\mu}}{\hat{\sigma}}$ goes to infinity for the positive side or equivalently negative side of these score functions in Eqs. \eqref{Scorea1a2a3a4nosign} or \eqref{Scorea1a2a3a4Huber} on the real line. In four designs, there are situations in which $\alpha$ is greater than $1$. Even if $\alpha_1,\alpha_3>1$, $\widehat{\text{MSE}}(\hat{\boldsymbol \theta})$, $\hat{\boldsymbol \theta}=(\hat{\mu},\hat{\sigma})$ for $\hat{\mu}$ and $\hat{\sigma}$ are small enough. In Tables \ref{sim2Design}-\ref{sim3Design}, values of $\widehat{\text{MSE}}(\hat{\sigma})$ for $\hat{\sigma}$ from $S^H$ as Huber's M-function do not decrease even if sample size $n$ is increased from $210$ to $410$, which can be an expected result, because Huber's M-function in Eq. \eqref{Huberpsi} depends $\alpha_2=2$ and $\alpha_1,\alpha_3=1$. For $2^{th}$ and $3^{rd}$ designs in subsection \ref{simdesigns}, we have $\alpha_2=3$ and $\alpha_1,\alpha_3=1.1,1.2,0.8$ which are different than Huber's M-function and it is not logical to expect that Huber's M-function fits other designs well if sample size is increased.  When we look at the simulation results for $1^{th}$ and $2^{nd}$ designs in Tables \ref{sim1Design}-\ref{sim2Design} respectively, these results show that main issue is about the modeling capability of score functions with EEs for the design of artificial data sets and the real data sets even if  $^\mu \psi_{\log}$ and $^\sigma \psi_{\log}$ for $\alpha>0$, $q=1$ are infinite in Eqs. \eqref{logqM} and \eqref{logqS} respectively (see also Table \ref{qvaluetable}). Since Eqs. \eqref{logqM}-\eqref{logqA} and \eqref{logDM}-\eqref{logDA} are similar expression except $w_i$, we can say same results for Eqs. \eqref{logDM}-\eqref{logDA} if $\beta=0$  (see also Table \ref{Dvaluetable}).  

The following results are observed generally in these designs: For the simultenaous estimations of parameters $\mu$ and $\sigma$,  MSE values of $\hat{\mu}$ and $\hat{\sigma}$ from score functions $S_q$ and $S^D$ with EEs have MSE values which are smaller than that of $S^H$, $S^{*H}$ and $S^{*}$ with EEs. MSE values of  $\hat{\mu}$ and $\hat{\sigma}$ from $S^{*H}$ and $S^{*}$ with EEs  have smaller than that of $S^H$ for sample sizes $210$ and $410$, because $S^H$ depends on $\alpha_2=2$ and $\alpha_1,\alpha_3=1$ for underlying and contamination, respectively. For simultenaous estimations of $\mu$, $\sigma$ and $\alpha$, MSE values of $\hat{\sigma}$ and $\hat{\alpha}$ from MDLE have smaller than that of MqLE. If we want to have small values of  $\widehat{\text{MSE}}(\hat{\sigma})$ and $\widehat{\text{MSE}}(\hat{\alpha})$ for $\hat{\sigma}$ and $\hat{\alpha}$, respectively, we need to change value of $q$ for each estimations of two parameters. Since parameter $q$ changes shape of function \cite{Bercher12a,Bercher12,Bercher13}, it is reasonable to observe this result.  $\beta$ in MDLE can be changed minimally to get small values of $\widehat{\text{MSE}}(\hat{\sigma})$ and $\widehat{\text{MSE}}(\hat{\alpha})$ for $\hat{\sigma}$ and $\hat{\alpha}$, respectively. For simultaneous estimation of $\mu$ and $\sigma$, the chosen values of TCs from $S^H$, $S^{*H}$ and $S^{*}$ can lead to have biased estimators for $\mu$, i.e.,   $\widehat{\text{Var}}(\hat{\mu})<\widehat{\text{MSE}}(\hat{\mu})$ in some cases from Tables \ref{sim1Design}-\ref{sim4Design}. For Tables \ref{sim1Design}-\ref{sim4Design}, bold represents the smallest values of $\widehat{\text{MSE}}(\hat{\boldsymbol \theta})$ among all of functions and italic represents the smallest values of $\widehat{\text{MSE}}(\hat{\boldsymbol \theta})$ among $S^H$, $S^{*H}$ and $S^{*}$.

The design of artificial data set affects success of simultenaous estimations for parameters of EP.  The different designs of artificial data sets were tried and they are not given in here due to the number of page restriction. In order to avoid the similarity of contamination schema in four designs applied at here, we chose arbitrary values of $\mu_1,\sigma_1$ and $\mu_3$ and $\sigma_3$ for contamination. We also tried same values of parameters for contamination. In this case, the values of TCs should be updated to have small values of $\widehat{\text{MSE}}(\hat{\boldsymbol \theta})$ for each tried cases of contamination. $k$, $t$, $q$ and $\beta$ should be chosen accurately for each sample sizes to have small MSE values when $n$ gets larger, because the designs in each sample sizes can be different than the other one. In other words, it is obvious that the values of TC change the structure of M-function. So the values of $\widehat{\text{MSE}}(\hat{\boldsymbol \theta})$ in different sample sizes are affected, as expected.

\subsection{Application to real data sets: Procedure of computation for estimation of $\mu, \sigma, \alpha$}\label{secreal}
The computation is performed according to the following order:
\begin{enumerate}
	\item    Arbitrary appropriate values of TCs and shape parameter $\alpha$ are chosen according to design of data set (outliers or contamination)  by user to start computation.
	\item    The values of $k$, $t$, $q$ and $\beta$ as TCs and shape parameters $\alpha_1$, $\alpha_2$ and $\alpha_3$ are chosen by grid search until the smallest value of volume in Eq. \eqref{entropygeneratorvol}, appropriate minimum values of IC in its self  in Eqs. \eqref{AICformula} - \eqref{AICformulalogq} and the smallest value of mean absolute error (MAE)  in Eq. \eqref{MAEeq} which is more precise evaluation for testing the fitting performance of objective and score functions are obtained.
	
	MAE is defined as 
	\begin{equation}\label{MAEeq}
	\frac{1}{n}  \sum_{i=1}^n |\text{ real}~ \text{sort}(\boldsymbol{x})  - \text{artificial}~ \text{sort}(\boldsymbol{x})|.
	\end{equation}
\end{enumerate}
After data are sorted, sample sizes of contamination and underlying distributions are determined.  $n_1=7, n_2=109$ and $ n_3=2$ for non-outlier case in real data sets. These sample sizes are determined by user according to outliers.
The estimates of  $\hat{\mu}$, $\hat{\sigma}$, $\hat{\alpha}$ for MqLE and MDLE, $\hat{\mu}$, $\hat{\sigma}$, fixed values of $\alpha_1,\alpha_2,\alpha_3$ for $S^H$, $S^{*H}$, $S^*$, $S^D$ and  $S_q$ are values used to generate random numbers from EP distribution.  Since redescending M-functions move slowly decreasing for tails when it is compared by  $S^H$, $S^{*H}$, $S^*$ \cite{Hampeletal86},  we use $\hat{\mu}$, $\hat{\sigma}$ and fixed values of $\alpha_1,\alpha_2,\alpha_3$ for $S_q$ and $S^D$.  Random numbers $x_j$ are obtained by $D_j$, as given by items for design of artificial data set in \ref{simdesigns}.   The replication number is $m=2 \cdot 10^4$ for $\sum_{j=1}^{3} n_j$. 

Since the main criteria for model selection among M-functions is MAE, we overcome the disadvantages occured by special function $\Gamma$ in FI matrices for volume and also tail part of $S$ and weighted score, i.e. $wS$ for IC. Tail part of $S$ and $wS$ can pull down IC hardly. For this reason, MAE  should be used to adjust appropriate values of IC. Volume and IC are tools to decrease high elapsed time occuring due to using MAE in several replications for computation. In other words, using only MAE consumes much time. Consulting volume and IC gives an advantage to decrease computation time while determining  values of TCs. Note that the elapsed time is around 1.5 hours which change according to sample size $n$.  
\section{Proof of FI matrices based on derivative of score functions  $S^H$, $S^*$,  $S^{*H}$, $S_q$ and $S^D$ from EEs }\label{proofofFisher}
The approach based on objective function $\rho$ \cite{AmariIG,KullbackIT,CanKor16} is adopted to score functions \cite{CanKor18,Plastinoetal97}. We can add or remove minus because relative entropy or divergence is used for the estimation of parameters $\boldsymbol{\theta}$. By using the definition of Tsallis $q$-entropy $H_q=\int_{a}^{b} f(x;\boldsymbol \theta)^q \log_q[f(x;\boldsymbol \theta)] \ud x$ for continuous variable $x$ and law of large numbers (LLN) \cite{LehmannCas98}, we can replace $f^q$ with $S$ and $\log_q$ with $S_q$ because $g(z;\boldsymbol{\theta})=S(x;\boldsymbol{\theta})^z$ in Eq. \eqref{partitonf} is taken. Further, Eq. \eqref{1qz} is equivalently rewritten as $f(x;\boldsymbol{\theta})\frac{f(x;\boldsymbol{\theta})^{1-q}-1}{1-q} = f(x;\boldsymbol{\theta}) \log_q(f(x;\boldsymbol{\theta})) \approxeq S(x;\boldsymbol{\theta})S_q(x;\boldsymbol{\theta})$ if Jackson $q$-derivative is ${g(qz;\boldsymbol{\theta})-g(z;\boldsymbol{\theta}) \over (q-1)z} \approxeq {g(z;\boldsymbol{\theta})-g(qz;\boldsymbol{\theta}) \over (1-q)z}  $. Similarly, we have $S$ and $S^D$  for log-likelihood and distorted log-likelihood, respectively.  Let us assume that score functions $S$, $S_q$ and $S^D$ are differentiable. Thus, entropy function can be based on score functions and FI matrix based on $S_q$ is proved by 
\begin{eqnarray}\label{DifKLPsiScore}
\widehat{F_{\log_q}}(S_q;f;\boldsymbol \theta)&=& \frac{1}{n} \sum_{i=1}^n \int_{(i-1)/n}^{i/n} \Delta S(x_i;\boldsymbol \theta) \Delta S_q(x_i;\boldsymbol \theta) \ud x_i, \\ \nonumber
&\stackrel{LLN}{=}& \sum_{i=1}^n \int_{(i-1)/n}^{i/n} f(x_i;\boldsymbol \theta) \Delta S(x_i;\boldsymbol \theta) \Delta S_q(x_i;\boldsymbol \theta) \ud x_i, \\ \nonumber
&=& \sum_{i=1}^n \int_{(i-1)/n}^{i/n}f(x_i;\boldsymbol \theta) \frac{\partial}{\partial \boldsymbol \theta}  S_q(x_i;\boldsymbol \theta)   \frac{\partial}{\partial \boldsymbol \theta} S(x_i;\boldsymbol \theta)^T \ud x_i, \\ \nonumber
&=& \sum_{i=1}^n \int_{(i-1)/n}^{i/n} f(x_i;\boldsymbol \theta)^{2-q}[(1-q)S^2(x_i;\boldsymbol \theta)+ \frac{\partial}{\partial \boldsymbol \theta} S(x_i;\boldsymbol \theta)] \frac{\partial}{\partial \boldsymbol \theta} S(x_i;\boldsymbol \theta)^T  \ud x_i.  \nonumber
\end{eqnarray}
Eq. \eqref{DifKLPsiScore} keeps spirit of Riemann integration or histograms as random bins on the real line. That is, let $\Delta(\boldsymbol{\theta})=\boldsymbol{\theta}-\hat{\boldsymbol{\theta}}=\boldsymbol{\theta}-(\boldsymbol{\theta}+h)$ be a difference operator. For $\lim_{h \rightarrow 0} \Delta(\boldsymbol{\theta}) =0$. $h \Delta(S)=S(x;\boldsymbol{\theta})-S(x;\hat{\boldsymbol{\theta}}) =S - \hat{S}$  as a derivative of $S$ w.r.t $\boldsymbol{\theta}$, i.e., $\frac{\partial}{\partial \boldsymbol{\theta}} S(x;\boldsymbol{\theta})=\lim_{h \rightarrow 0} \Delta(S)$. FI based on the distorted log-score function $S^D$ is obtained similarly. For all intervals $x_i$ of sampled version of the fourth line of Eq. \eqref{DifKLPsiScore}, we have result in Eq. \eqref{derivSFisherq} for continuous case of $x$ with interval $(a,b)$.

\subsection{Proof for positive semidefinite of Fisher information matrices in Eqs. \eqref{derivSFisher}, \eqref{derivSFisherq} and \eqref{derivSFisherbeta}}\label{proofofFisherM}

If we produce EEs of location $\mu$ and scale $\sigma$ as a location-scale family in a p.d. function $f(x;\boldsymbol{\theta})$, score function $S$ and $\Psi$ for $\mu$ and $\sigma$ can have same mathematical expression, as given by Eqs. \eqref{scorelogf} and \eqref{Huberpsi}. For location-scale model, $\Psi$ for $\mu$ and $\sigma$ can drop to $S$ as only one function for $\mu$ and $\sigma$ due to the spirit of EEs, as introduced by subsection \ref{introEEs}. Since we can have EE for shape parameter $\alpha$ of EP distribution, only one score function for $\mu$, $\sigma$ and $\alpha$ is obtained from EEs. Since $\rho(x;\boldsymbol{\theta})=-\log(f(x;\boldsymbol{\theta}))$,  $S(x;\boldsymbol{\theta})$ from $\Psi(x;\boldsymbol{\theta})$  is mainly from $\frac{f^{'}(x;\boldsymbol{\theta})}{f(x;\boldsymbol{\theta})}$. $f^{'}(x;\boldsymbol{\theta})=\frac{\partial}{\partial \boldsymbol{\theta}}f(x;\boldsymbol{\theta})$,  $f^{''}(x;\boldsymbol{\theta})=\frac{\partial^2}{\partial \boldsymbol{\theta}\partial \boldsymbol{\theta}^T}f(x;\boldsymbol{\theta})$, $\boldsymbol{\theta} \in (\mu,\sigma)$ \cite{Hub81}. Let $f(x;\boldsymbol{\theta})$, $f^{'}(x;\boldsymbol{\theta})$ and $f^{''}(x;\boldsymbol{\theta})$ be shown as $f$, $f^{'}$ and $f^{''}$, respectively. Since $S$ is rooted from $f^{'}/f$, diagonal elements of FI matrix in Eq. \eqref{derivSFisher} without underlying $f_2 \geq 0$ distribution is rewritten as following form:

\begin{equation}\label{scoresquareS2}
\left[\left(\frac{f^{'}}{f}\right)^{'}\right]^2   = \left[\frac{f^{''}}{f}-\left(\frac{f^{'}}{f}\right)^2\right]^2. 
\end{equation}
Since Eq. \eqref{scoresquareS2} has a square form, we have positive values for the diagonal elements of Fisher matrix $F_{\log}$ in  Eq. \eqref{derivSFisher}.    Let $f_{\mu}^{'}$ and  $f_{\mu}^{''}$ be $\frac{\partial}{\partial \mu}f(x;\mu)$ and $\frac{\partial^2}{\partial \mu^2}f(x;\mu)$, respectively. $f_{\mu}$ given by $\frac{\partial}{\partial \mu}\log(f(x;\mu,\sigma,\alpha))=\frac{f_{\mu}^{'}}{f_{\mu}}$ is p.d. function in $\log$. Similarly, other twin parameters, i.e., $\sigma$ and $\alpha$, etc. can be chosen to have notations in Eq. \eqref{scoresquareSS}.  The undiagonal elements of FI matrix in Eq. \eqref{derivSFisher}  without underlying $f_2 \geq 0$ distribution can be negative and positive due to the given following form:
\begin{equation}\label{scoresquareSS}
\left(\frac{f_{\mu}^{'}}{f_{\mu}}\right)^{'}\left(\frac{f_{\sigma}^{'}}{f_{\sigma}}\right)^{'}=\left[\frac{f_{\mu}^{''}}{f_{\mu}}-\left(\frac{f_{\mu}^{'}}{f_{\mu}}\right)^2 \right]  \left[\frac{f_{\sigma}^{''}}{f_{\sigma}}-\left(\frac{f_{\sigma}^{'}}{f_{\sigma}}\right)^2\right]. 
\end{equation}

The following conditions should be satisfied in order to imply that $F_{\log}$ in Eq. \eqref{derivSFisher} is positive semidefinite. Let us test $F_{\log}$ by use of these conditions. 

\begin{itemize}
	\item Determinants test:
	\begin{itemize}
		\item 	$E_{\mu\mu}=E\left((\frac{\partial}{\partial \mu}S(x;\mu,\sigma))^2\right) \geq 0$ and $E_{\sigma\sigma}=E\left((\frac{\partial}{\partial \sigma}S(x;\mu,\sigma))^2\right) \geq 0$ because of Eq \eqref{scoresquareS2}. 
		\item  $E_{\mu\mu}E_{\sigma\sigma} \geq E_{\mu\sigma}^2$, i.e. $E\left((\frac{\partial}{\partial \mu}S(x;\mu,\sigma))^2\right)E\left((\frac{\partial}{\partial \sigma}S(x;\mu,\sigma))^2\right) \geq \left(E\frac{\partial}{\partial \mu}S(x;\mu,\sigma)\frac{\partial}{\partial \sigma}S(x;\mu,\sigma)\right)^2$.
	\end{itemize} 
	\item Pivot test: $E_{\mu\mu}=E\left((\frac{\partial}{\partial \mu}S(x;\mu,\sigma))^2\right) \geq 0$ and $E_{\mu\mu}^{-1}(E_{\mu\mu}E_{\sigma\sigma}-E_{\mu \sigma}^2) \geq 0$.
\end{itemize}
Consequently, $F_{\log}$ in Eq. \eqref{derivSFisher} is  a positive semidefinite matrix.

Let us examine whether or not diagonal elements of FI matrices in Eqs. \eqref{derivSFisherq} and Eq. \eqref{derivSFisherbeta} without $f^{2-q} \geq 0$ and $\frac{f^2}{\beta+f} \geq 0$ are positive.  Since $wS$ is rooted from $wf^{'}/f$, diagonal elements are rewritten as following form:
\begin{equation}\label{scoresquareS3}
\frac{(f^{'})^4}{f^{4}} q - f^{''}\frac{(f^{'})^2}{f^{3}}(q+1)+\frac{(f^{''})^2}{f^{2}}.
\end{equation}
\begin{equation}\label{scoresquareS4}
\frac{(f^{'})^4}{f^{3}} \frac{1}{\beta+f} - f^{''}\frac{(f^{'})^2}{f^3}\left(2-\frac{\beta}{\beta+f}\right)+\frac{(f^{''})^2}{f^2},
\end{equation}
\noindent $\beta \geq 0$ and $q \in (0,1)$. If $f$ is differentiable and concave according to parameters, then $f^{''}$ is negative. Thus, Eqs. \eqref{scoresquareS3} and \eqref{scoresquareS4} are positive because other expressions in Eqs. \eqref{scoresquareS3} and \eqref{scoresquareS4} are positive or if summation of first and third terms in Eq. \eqref{scoresquareS3} is bigger than second term in Eq. \eqref{scoresquareS3}, then Eq. \eqref{scoresquareS3} is positive. The same rule is valid for Eq. \eqref{scoresquareS4}.  Note that since $f$ is p.d. function and $\beta \geq 0$, then $\beta+f \geq 0$ and $0 \leq \frac{\beta}{\beta+f} \leq 1$.

Let us examine undiagonal elements of FI matrices in Eqs. \eqref{derivSFisherq} and Eq. \eqref{derivSFisherbeta} without $f^{2-q} \geq 0$ and $\frac{f^2}{\beta+f} \geq 0$, as given by following forms, respectively.
\begin{equation}\label{scoresquareS5}
\left(\frac{f_{\sigma}^{''}}{f_{\sigma}}-\left(\frac{f_{\sigma}^{'}}{f_{\sigma}}\right)^2\right) \left((1-q)\left(\frac{f^{'}}{f}\right)^2+\frac{f_{\mu}^{''}}{f_{\mu}} - \left(\frac{f_{\mu}^{'}}{f_{\mu}}\right)^2\right).
\end{equation}
\begin{equation}\label{scoresquareS6}
\left(\frac{f_{\sigma}^{''}}{f_{\sigma}}-\left(\frac{f_{\sigma}^{'}}{f_{\sigma}}\right)^2\right) \left(\left(\frac{\beta}{\beta+f}\right)\left(\frac{f^{'}}{f}\right)^2+\frac{f_{\mu}^{''}}{f_{\mu}} - \left(\frac{f_{\mu}^{'}}{f_{\mu}}\right)^2\right).
\end{equation}
Eqs. \eqref{scoresquareS5} and \eqref{scoresquareS6} can be positive or negative according to values of expressions. $\mu$ and $\sigma$  can be replaced by $\sigma$ and $\alpha$, respectively. Other possible replacement for parameters can be done to have undiagonal elements of Fisher matrix in Eq. \eqref{musigmaalphaF}. Thus,  undiagonal elements of inverse of FI matrices in Eqs. \eqref{derivSFisherq} and \eqref{derivSFisherbeta} can be positive or negative.  Note that $0 \leq f_2 \leq 1 $, $0 \leq f^{2-q} \leq 1$ and $\frac{f^2}{\beta+f} \geq 0$ do not change the sign of expectation or expressions in Eqs. \eqref{scoresquareS2}-\eqref{scoresquareS6}.

$F_{\log_q}$ and $F_{\log}$ in Eqs. \eqref{derivSFisherq} and \eqref{derivSFisherbeta} will be positive semidefinite matrices if square of Eq. \eqref{scoresquareS3} is greater and equal than multiplication of Eq. \eqref{scoresquareS5} and Eq. \eqref{scoresquareS5} with replacement of $\sigma$ and $\mu$ according to matrix in Eq. \eqref{musigmaalphaF}, i.e., determinants of upper submatrices are positive or zero if there exists equality. The same procedure given for Eqs. \eqref{scoresquareS3} and \eqref{scoresquareS5} is valid for \eqref{scoresquareS4} and \eqref{scoresquareS6}. $F_{\log_q}$ and $F_{\log}$ are not symmetric, but $F$ is a square matrix with $d \times d$. If Eqs. \eqref{scoresquareS3} and \eqref{scoresquareS4} are positive, then diagonal elements of $F$ will be positive. Thus, we have positive value for trace of $F$, i.e. summation of eigenvalues of $F$ is positive if the number of columns or rows is equal to rank of $F$. Determinant of $F$ is product of eigenvalues of $F$.  When every eigenvalue of $F$ is positive, $F$ is positive definite due to eigenvalue decomposition if $F$ has $d$ linearly independent eigenvectors \cite{matrixproperty}.  If determinant of $F$ is zero for the positive semidefinite case, then the Moore-Penrose generalized inverse for singular matrix  \cite{LiYeh12gen} can be used to get the variances of estimators, i.e. Var$(\hat{\mu})$, Var$(\hat{\sigma})$ and Var$(\hat{\alpha})$, even if information loss occurring due to  using the generalized inverse exists. 

To show whether or not $F$ is positive semidefinite, we have used $\frac{\partial}{\partial \boldsymbol{\theta}}\rho(x;\boldsymbol{\theta})=\frac{\partial}{\partial \boldsymbol{\theta}} \Lambda \cdot f(x;\boldsymbol{\theta})^{'}=\frac{\partial}{\partial \boldsymbol{\theta}} \Lambda \cdot  {f(x;\boldsymbol{\theta})^{'} \over f(x;\boldsymbol{\theta})} f(x;\boldsymbol{\theta})=w {f(x;\boldsymbol{\theta})^{'} \over f(x;\boldsymbol{\theta})}$. The equivalence between $wS$ and $\frac{\partial}{\partial \boldsymbol{\theta}}\rho=\Psi$ is enjoyed in order to be able to say that we have  positive semidefinite $F$ tentatively for $wS$ used in FI. Otherwise, as we already know that estimation process is based on divergence which is mainly  an absolute value of distance between $f(x;\boldsymbol{\theta})$ and $f(x;\hat{\boldsymbol{\theta}})$ for p.d. function \cite{PardoSD,KullbackIT} or equivalently $\rho(x;\boldsymbol{\theta})$ and $\rho(x;\hat{\boldsymbol{\theta}})$  for objective functions or equivalently  $S(x;\boldsymbol{\theta})$ and $S(x;\hat{\boldsymbol{\theta}})$ for score functions, elements of FI matrix can be multiplied by minus accordingly in order to avoid negative values for diagonal elements of inverse of FI matrix in which we can have for score functions used only in FI.

\section{Fisher information matrices for the parameters $\mu,\sigma$ and $\alpha$ of EP distribution}
FI matrices based on score functions in Eq. \eqref{derivSFisher} for $\mu$ and $\sigma$ and also Eqs. \eqref{derivSFisherq}-\eqref{derivSFisherbeta} for $\mu$, $\sigma$ and $\alpha$ are represented by the following forms with elements:
\begin{equation}\label{musigmaalphaF}
F=n\begin{bmatrix} 
E_{\mu \mu} & E_{\mu \sigma} & E_{\mu \alpha} \\
E_{\sigma \mu} & E_{\sigma \sigma} & E_{\sigma \alpha}  \\
E_{\alpha \mu} & E_{\alpha \sigma} & E_{\alpha \alpha}  
\end{bmatrix},
\end{equation}
\noindent where $E$ represents an expectation taken over the underlying distribution for derivative of score functions w.r.t parameters $\mu$, $\sigma$ and $\alpha$. $n$ is sample size \cite{Hub64}.
\subsection{The elements of FI matrix in Eq. \eqref{derivSFisher} based on derivative of combined $\log$-score $S^{*}$ from EEs for parameters of EP distribution if $f_2$ is an underlying distribution}
\begin{equation}\label{musigmaalphaFS}
F_{\log}=n\begin{bmatrix} 
E_{\mu \mu}\left[\left(\frac{\partial S}{\partial \mu}\right)^2 \right] & E_{\mu \sigma}\left[\frac{\partial S}{\partial \mu}\frac{\partial S}{\partial \sigma}\right]  \\
E_{\sigma \mu}\left[\frac{\partial S}{\partial \sigma}\frac{\partial S}{\partial \mu}\right] & E_{\sigma \sigma}\left[\left(\frac{\partial S}{\partial \sigma}\right)^2 \right]   
\end{bmatrix},
\end{equation}
\noindent where $E_{\mu \sigma}\left[\frac{\partial S}{\partial \mu}\frac{\partial S}{\partial \sigma}\right]=E_{\sigma \mu}\left[\frac{\partial S}{\partial \sigma}\frac{\partial S}{\partial \mu}\right]$. Eqs. \eqref{Esmsm} and \eqref{Esgsg} are positive. As shown by Eq.\eqref{scoresquareSS}, Eq. \eqref{Esmsg} can be negative and positive. Thus, matrix in Eq. \eqref{musigmaalphaFS} is positive semidefinite due to tests of determinants and pivot.

The reparametrization of $\Gamma$ function  is used to calculate the integrals in Eqs. \eqref{Esmsm}-\eqref{Esgsg} \cite{CanEnt18}. Note that $k$ and $t$ should be positive due to $\Gamma$ \cite{Gradshteyn07}.
\begin{footnotesize}
\begin{eqnarray}\label{Esmsm}
E_{\mu \mu}\left[\left(\frac{\partial S}{\partial \mu}\right)^2 \right]&=&\frac{1}{2\sigma^2\Gamma(1/\alpha_2)}\{(\alpha_1^2-\alpha_1)^2\Gamma(\frac{2\alpha_1-3}{\alpha_2},(k/\sigma)^{\alpha_2})+(\alpha_3^2-\alpha_3)^2\Gamma(\frac{2\alpha_3-3}{\alpha_2},(t/\sigma)^{\alpha_2}) \\ \nonumber
&+&(\alpha_2^2-\alpha_2)^2\left[\gamma(2-3/\alpha_2,(k/\sigma)^{\alpha_2})+\gamma(2-3/\alpha_2,(t/\sigma)^{\alpha_2})\right]\},
\end{eqnarray}
\end{footnotesize}
\begin{footnotesize}
\begin{eqnarray}\label{Esmsg}
E_{\mu \sigma}\left[\frac{\partial S}{\partial \mu}\frac{\partial S}{\partial \sigma}\right]&=&\frac{1}{2\sigma^2\Gamma(1/\alpha_2)}\{-(\alpha_1^2-\alpha_1)^2\Gamma(\frac{2\alpha_1-2}{\alpha_2},(k/\sigma)^{\alpha_2})+(\alpha_3^2-\alpha_3)^2\Gamma(\frac{2\alpha_3-2}{\alpha_2},(t/\sigma)^{\alpha_2}) \\ \nonumber
&+&(\alpha_2^2-\alpha_2)^2\left[-\gamma(2-3/\alpha_2,(k/\sigma)^{\alpha_2})+\gamma(2-3/\alpha_2,(t/\sigma)^{\alpha_2})\right]\},
\end{eqnarray}
\end{footnotesize}
\begin{footnotesize}
\begin{eqnarray}\label{Esgsg}
E_{\sigma \sigma}\left[\left(\frac{\partial S}{\partial \sigma}\right)^2 \right]&=&\frac{1}{2\sigma^2\Gamma(1/\alpha_2)}\{(\alpha_1^2-\alpha_1)^2\Gamma(\frac{2\alpha_1-1}{\alpha_2},(k/\sigma)^{\alpha_2})+(\alpha_3^2-\alpha_3)^2\Gamma(\frac{2\alpha_3-1}{\alpha_2},(t/\sigma)^{\alpha_2}) \\ \nonumber
&+&(\alpha_2^2-\alpha_2)^2\left[\gamma(2-1/\alpha_2,(k/\sigma)^{\alpha_2})+\gamma(2-1/\alpha_2,(t/\sigma)^{\alpha_2})\right]\}.
\end{eqnarray}
\end{footnotesize}
$ \Gamma(z)=\gamma(z,a) + \Gamma(z,a)$ is gamma function with lower and upper incomplete gamma functions \cite{Gradshteyn07}. For $S^{*H}$, $S_1$ and $S_3$ in $S^*$ are multiplied by $-k$ and $t$, respectively. $\alpha_1,\alpha_2,\alpha_3>3/2$ due to $\Gamma$ function. 

\subsection{The elements of FI matrix in Eq. \eqref{derivSFisherq} based on derivative of score $S_q$ from EEs for parameters of EP distribution if $f$ is an underlying distribution}

\begin{equation}\label{musigmaalphaFlogq}
F_{\log_q}=n\begin{bmatrix} 
E_{\mu \mu}\left[(1-q)S^2 \frac{\partial S}{\partial \mu}+\frac{\partial S}{\partial \mu}\frac{\partial S}{\partial \mu} \right] & E_{\mu \sigma}\left[ (1-q)S^2 \frac{\partial S}{\partial \sigma}+\frac{\partial S}{\partial \mu}\frac{\partial S}{\partial \sigma}\right] & E_{\mu \alpha}\left[ (1-q)S^2 \frac{\partial S}{\partial \alpha}+\frac{\partial S}{\partial \mu}\frac{\partial S}{\partial \alpha}\right] \\
E_{\sigma \mu}\left[ (1-q)S^2 \frac{\partial S}{\partial \mu}+\frac{\partial S}{\partial \sigma}\frac{\partial S}{\partial \mu}\right] & E_{\sigma \sigma}\left[ (1-q)S^2 \frac{\partial S}{\partial \sigma}+\frac{\partial S}{\partial \sigma}\frac{\partial S}{\partial \sigma}\right] & E_{\sigma \alpha}\left[ (1-q)S^2 \frac{\partial S}{\partial \alpha}+\frac{\partial S}{\partial \sigma}\frac{\partial S}{\partial \alpha}\right]  \\
E_{\alpha \mu}\left[ (1-q)S^2 \frac{\partial S}{\partial \mu}+\frac{\partial S}{\partial \alpha}\frac{\partial S}{\partial \mu}\right] & E_{\alpha \sigma}\left[ (1-q)S^2 \frac{\partial S}{\partial \sigma}+\frac{\partial S}{\partial \alpha}\frac{\partial S}{\partial \sigma}\right] & E_{\alpha \alpha}\left[ (1-q)S^2 \frac{\partial S}{\partial \alpha}+\frac{\partial S}{\partial \alpha}\frac{\partial S}{\partial \alpha}\right]  
\end{bmatrix},
\end{equation}
The matrix $F_{\log_q}$  in Eq. \eqref{musigmaalphaFlogq} and its form with parameters $\mu$ and $\sigma$ as two-dimensional matrix are positive semidefinite if conditions are satisfied  from Eqs. \eqref{scoresquareS3} and \eqref{scoresquareS5}. 

Mathematica 11.3 is used to calculate the integrals in Eqs. \eqref{Esmsmq}-\eqref{Esalalq}.
\begin{small}
\begin{eqnarray} \label{Esmsmq}
E_{\mu \mu}\left[(1-q)S^2 \frac{\partial S}{\partial \mu}+\frac{\partial S}{\partial \mu}\frac{\partial S}{\partial \mu} \right]&=&-2^{q-1}(\alpha-1)\alpha^{3-q}\sigma^{q-3}(2-q)^{3/\alpha-3} \\ \nonumber
&&\left(2+3\sigma(q-1)-q+\alpha(2+2\sigma(1-q)+q)\right) \Gamma(2-3/\alpha)\Gamma(1/\alpha)^{q-2}
\end{eqnarray}
\end{small}
\begin{eqnarray}\label{Esmsgq}
E_{\mu \sigma}\left[ (1-q)S^2 \frac{\partial S}{\partial \sigma}+\frac{\partial S}{\partial \mu}\frac{\partial S}{\partial \sigma}\right]&=&0,
\end{eqnarray}
\begin{eqnarray}\label{Esmalq}
E_{\mu \alpha}\left[ (1-q)S^2 \frac{\partial S}{\partial \alpha}+\frac{\partial S}{\partial \mu}\frac{\partial S}{\partial \alpha}\right]&=&0,
\end{eqnarray}
\begin{eqnarray}\label{Esgsmq}
E_{\sigma \mu}\left[ (1-q)S^2 \frac{\partial S}{\partial \mu}+\frac{\partial S}{\partial \sigma}\frac{\partial S}{\partial \mu}\right]&=&\frac{\alpha^{4-q}(1-q)(1-\alpha)\Gamma(3-3/\alpha)}{2^{1-q}(\Gamma(1/\alpha)\sigma)^{2-q}(2-q)^{3-3/\alpha}},
\end{eqnarray}
\begin{eqnarray}\label{Esgsgq}
E_{\sigma \sigma}\left[ (1-q)S^2 \frac{\partial S}{\partial \sigma}+\frac{\partial S}{\partial \sigma}\frac{\partial S}{\partial \sigma}\right]&=&
2^{q-2}(\alpha-1)^2\alpha^{1-q}\sigma^{q-3}(2-q)^{1/\alpha-2} \\ \nonumber
&&\left(\alpha^2\Gamma(2-1/\alpha)-(\alpha-1)\Gamma(-1/\alpha)\right)\Gamma(1/\alpha)^{q-2}
\end{eqnarray}
\begin{eqnarray}\label{Essgalq}
E_{\sigma \alpha}\left[ (1-q)S^2 \frac{\partial S}{\partial \alpha}+\frac{\partial S}{\partial \sigma}\frac{\partial S}{\partial \alpha}\right]&=&(2-q)^{1/\alpha-2}2^{q-1}(\alpha-1)\Gamma(2-1/\alpha)\left(\frac{\alpha}{\sigma\Gamma(1/\alpha)}\right)^{2-q} \\ \nonumber
&&\left(\log(2-q)-\psi^{(0)}(2-1/\alpha)-1\right),
\end{eqnarray}
\begin{tiny}
\begin{eqnarray}\label{Esalmq}
E_{\alpha \mu}\left[ (1-q)S^2 \frac{\partial S}{\partial \mu}+\frac{\partial S}{\partial \alpha}\frac{\partial S}{\partial \mu}\right]&=&2^{q-1}(\alpha-1)\alpha^{4-q}\sigma^{q-2}(2-q)^{3/\alpha-3}(q-1)\Gamma(3-3/\alpha)\Gamma(1/\alpha)^{q-2}
\end{eqnarray}
\end{tiny}
\begin{small}
\begin{eqnarray}\label{Esalsgq}
E_{\alpha \sigma}\left[ (1-q)S^2 \frac{\partial S}{\partial \sigma}+\frac{\partial S}{\partial \alpha}\frac{\partial S}{\partial \sigma}\right]&=& 2^{q-1}(\alpha-1)(\sigma/\alpha)^{q-2}(2-q)^{1/\alpha-2}\Gamma(2-1/\alpha)\Gamma(1/\alpha)^{q-2} \\ \nonumber 
&&\left(\log(2-q)-1+\psi^{(0)}(2-1/\alpha)\right)
\end{eqnarray}
\end{small}
\begin{tiny}
\begin{eqnarray}\label{Esalalq}
E_{\alpha \alpha}\left[ (1-q)S^2 \frac{\partial S}{\partial \alpha}+\frac{\partial S}{\partial \alpha}\frac{\partial S}{\partial \alpha}\right]&=&
(2\sigma/\alpha)^{q-1}(2-q)^{1/\alpha-2}\Gamma(2-1/\alpha)\Gamma(1/\alpha)^{q-2} \\ \nonumber
&&\left((\log(2-q)-1)^2-2(\log(2-q)-1)\psi^{(0)}(2-1/\alpha)+\psi^{(0)}(2-1/\alpha)^2+\psi^{(1)}(2-1/\alpha)\right)
\end{eqnarray}
\end{tiny}
For $q=1$, FI based on $S$ is obtained by Eqs. \eqref{Esmsmq}-\eqref{Esalalq}. $\psi$ is digamma function and $\psi^{(h)}$ is $h^{th}$ derivative of the digamma function. $\alpha>3/2$ is taken because of $\Gamma$ function and $q \in (0,1)$. 

The numerical integration in MATLAB R2013a is used for elements of FI matrix in Eq.  \eqref{derivSFisherbeta}. Since Eq. \eqref{derivSFisherbeta} is same with Eq. \eqref{derivSFisherq}, the matrix in Eq. \eqref{derivSFisherbeta} or Eq. \eqref{musigmaalphaFlogq} based on $S^D$ is positive semidefinite if conditions (positive and negative{}) for Eqs. \eqref{scoresquareS4} and \eqref{scoresquareS6} are satisfied. For parameters $\mu$ and $\sigma$, we have  two-dimensional matrix of Eq. \eqref{musigmaalphaFlogq}.

\section*{Acknowledgements}
We would like to thank so much Editorial Board and anonymous referees to provide the invaluable comments.  
\section*{Disclosure statement}
No potential conflict of interest was reported by the author(s).

\end{document}

%% file: table1.tex
\begin{table}[!htb] 
	\centering
	\caption{Limit values of $\Psi_{\log_q}		=( \underset{x \rightarrow \infty}{\lim} ^\mu \psi_{\log_q},  \underset{x \rightarrow \infty}{\lim} ^\sigma \psi_{\log_q},   \underset{x \rightarrow \infty}{\lim} ^\alpha \psi_{\log_q}) $.}
	\begin{tabular}{c|cccc}
		& $q=1$ & $0<q<1$ & $q>1$ &  \\ \hline
		$\alpha=1$	& $(-\alpha,\infty,-\infty)$ & $(0,0,0)$ & $(-\infty,\infty,-\infty)$ &  \\
		$0<\alpha<1$  & $(0,\infty,-\infty)$  & $(0,0,0)$ & $(\infty,\infty,-\infty)$  \\
		$\alpha>1$& $(-\infty,\infty,-\infty)$ &  $(0,0,0)$ &  $(-\infty,\infty,-\infty)$  & 
	\end{tabular}
	\label{qvaluetable}
\end{table}

%% file: table2.tex
\begin{table}[!htb]\centering
	\caption{Limit values of $\Psi_{\log}^D=(	\underset{x \rightarrow \infty}{\lim} ^\mu \psi_{\log}^D, \underset{x \rightarrow \infty}{\lim} ^\sigma \psi_{\log}^D, \underset{x \rightarrow \infty}{\lim} ^\alpha \psi_{\log}^D   )$.}
	\begin{tabular}{c|cccc}
		& $\alpha=1$ & $0<\alpha<1$ & $\alpha>1$ &  \\ \hline
		$\beta>0$ & $(0,0,0)$ & $(0,0,0)$ & $(0,0,0)$ &  \\
		$\beta=0$ & $(-\alpha,\infty,-\infty)$  & $(0,\infty,-\infty)$ & $(-\infty,\infty,-\infty)$  
	\end{tabular}
	\label{Dvaluetable} 
\end{table}

%% file: table3.tex
\begin{table}[!htb]
	\centering
	\caption{Example 1: M-estimates of parameters $\boldsymbol \theta=(\mu,\sigma,\alpha)$ }
	\label{ex1temperature}
	\scalebox{0.64}
	{
		\begin{tabular}{ccccccccccccc}\hline
			EEs & TC & $\hat{\mu}$ & $\hat{\sigma}$ &$\alpha_1,\alpha_2,\alpha_3$  & $\text{Var}_S(\hat{\mu})$ & $\text{Var}_S(\hat{\sigma})$  & $\text{Var}_S(\hat{\alpha})$  & $\text{AIC}_S$ & $\text{cAIC}_S$ & $\text{BIC}_S$ & $\text{Vol}(S;\hat{\mu},\hat{\sigma},\hat{\alpha})$ & MAE   \\ \hline 			
			$S^H$ & $k=.95$, $t=4.1$ & 4.2000 & 1.5108 & 1,2,1 &.0148 &.0141 & - &  186.5247 &  186.6552 & 191.6325 &.001896&1.0334  \\  
			outliers&   & 4.2000& 1.6172 & 1,2,1 &.0164 &.0158 & - &  196.2415 & 196.3692 & 201.3909&.002060&1.2010 \\	 \\ 
			$S^{*H}$&$k=0$, $t=3.15$ & 3.6135 &2.3078 & 1.25,2.1,1.4 &.0307&.0261 & - & 52.8691 &  52.9996 &  57.9769 & .003497&0.5214   \\
			outliers & & 3.7006 & 2.7360 & 1.25,2.1,1.4 &.0357  &.0313 & - &55.2169 &  55.3445  & 60.3663 & .004179&0.8245  \\	 \\ 
			$S^{*}$&$k=.95$, $t=4$ & 3.0674 & 1.4233 & 1.25,2.1,1.4  &.0281&.0258 & - & 253.9906 & 254.1210  & 259.0983 & .002632&0.1912   \\
			outliers & &  3.1192 &  1.8391 & 1.25,2.1,1.4 &.0462 &.0452 & - &225.5623 & 225.6899&  230.7117& .004141&$\boldsymbol{0.3157}$  \\	 \hline
			EEs, $\rho$ & TC & $\hat{\mu}$ & $\hat{\sigma}$ &$\alpha, \hat{\alpha}$  & $\text{Var}_S(\hat{\mu})$ & $\text{Var}_S(\hat{\sigma})$  & $\text{Var}_S(\hat{\alpha})$  & $\text{AIC}_S$ & $\text{cAIC}_S$ & $\text{BIC}_S$ & $\text{Vol}(S;\hat{\mu},\hat{\sigma},\hat{\alpha})$ & MAE   \\ \hline 		
			$S^{D}$&$\beta=10^{-2}$ & 3.1201& 1.6752& 2.1 &.0089&.0139& - &197.8566 & 197.9870 & 202.9643 & .001475&0.1162   \\
			outliers & &  3.1201 &1.6752& 2.1  &.0087 &.0136& - &197.8566 & 197.9842&  203.0060 & .001414&0.3269  \\	\\
			MDLE&$\beta=10^{-2}$ & 3.1251 & 1.6264 & 1.9026  &.0095&.0394 & .0270 &198.3792 & 198.6429 & 206.0408 &.000155&$\boldsymbol{0.1055}$   \\
			outliers & &  3.1261 &  1.6374 & 1.9375  &.0093 &.0361 & .0272 & 198.5424 & 198.8005 & 206.2666 & .000142&0.3337 \\	  \hline
			EEs, $\rho$ & TC & $\hat{\mu}$ & $\hat{\sigma}$  &$\alpha, \hat{\alpha}$ & $\text{Var}_{S_q}(\hat{\mu})$ & $\text{Var}_{S_q}(\hat{\sigma})$ & $\text{Var}_{S_q}(\hat{\alpha})$ & $\text{AIC}_{S_q}$ & $\text{cAIC}_{S_q}$ &  $\text{BIC}_{S_q}$ & $\text{Vol}_q(S_q;\hat{\mu},\hat{\sigma}, \hat{\alpha})$ & MAE   \\ \hline 
			$S_q$& $q=.8$ & 3.1358& 1.6018 & 2.1 &.0117 &.0173 & -  & 153.6336& 153.7640 & 158.7413  & .001876&0.1231   \\
			outliers &  &3.1358 &1.6018 & 2.1  &.0114 &.0169 &- &153.6355&  153.7632 & 158.7849 & .001799&0.3431 \\ \\	 
			MqLE& $q=.8$ & 3.1345&1.5386&1.9054 &.0116 &.0440& .0336 & 155.0965 & 155.3602 & 162.7581&.000211&0.1238 \\
			outliers &  & 3.1385 &1.5404 & 1.9311 &.0114 &.0407 & .0336 &155.7167 & 155.9748 & 163.4408&.000193&0.3657 \\	\hline
			
		\end{tabular}
	}
\end{table}

%% file: table4.tex
\begin{table}[!htb]
	\centering
	\caption{Example 2: M-estimates of parameters $\boldsymbol \theta=(\mu,\sigma,\alpha)$}
	\label{ex2microarray}
	\scalebox{0.6}
	{
		\begin{tabular}{ccccccccccccc}\hline
			EEs & TC & $\hat{\mu}$ & $\hat{\sigma}$ &$\alpha_1,\alpha_2,\alpha_3$  & $\text{Var}_S(\hat{\mu})$ & $\text{Var}_S(\hat{\sigma})$  & $\text{Var}_S(\hat{\alpha})$  & $\text{AIC}_S$ & $\text{cAIC}_S$ & $\text{BIC}_S$ & $\text{Vol}(S;\hat{\mu},\hat{\sigma},\hat{\alpha})$ & MAE   \\ \hline 			
			$S^H$ & $k=-1.22$, $t=1.18$ & 0.0307 & 0.0677 & 1,2,1 & $0.9702\cdot10^{-4}$  &$0.9702\cdot10^{-4}$ & - & 234.8350 & 234.9393 & 240.3764 &$ 1.0308\cdot10^{-5}$&0.0869  \\  
			outliers&   & 0.0263& 0.1103 & 1,2,1 &$0.2533\cdot10^{-3}$ &$0.2533\cdot10^{-3}$ & - &204.0155 & 204.1181 & 209.5905&$  2.6457\cdot10^{-5}$&0.0746 \\	 \\ 
			$S^{*H}$&$k=-0.69$, $t=0.67$ & -0.0304 &0.1148 & 1.52,3.18,1.11 &$ 0.7908\cdot10^{-3}$&$ 0.7908\cdot10^{-3}$ & - & 175.0141 & 175.1184 & 180.5555 & $ 8.4022\cdot10^{-5}$&0.0669   \\
			outliers & & -0.0275 & 0.1379 & 1.52,3.18,1.11 &0.0011  &0.0011 & - &170.9423 & 171.0448&  176.5172 &  $ 1.1723\cdot10^{-4}$ &0.0695  \\	 \\ 
			$S^{*}$&$k=-0.65$, $t=0.61$ & -0.0255 & 0.1504 & 1.52,3.18,1.11  &$0.9220\cdot10^{-3}$&$0.9220\cdot10^{-3}$ & - & 208.4176&  208.5219 & 213.9589 & $ 9.8055 \cdot10^{-5}$&0.0488   \\
			outliers & &  -0.0223 &  0.1771 & 1.52,3.18,1.11 &0.0013 &0.0013 & - &  196.5984 & 196.7010 & 202.1734 & $ 1.3147\cdot10^{-4}$&0.0514  \\	 \hline
			EEs, $\rho$ & TC & $\hat{\mu}$ & $\hat{\sigma}$ &$\alpha, \hat{\alpha}$  & $\text{Var}_S(\hat{\mu})$ & $\text{Var}_S(\hat{\sigma})$  & $\text{Var}_S(\hat{\alpha})$  & $\text{AIC}_S$ & $\text{cAIC}_S$ & $\text{BIC}_S$ & $\text{Vol}(S;\hat{\mu},\hat{\sigma},\hat{\alpha})$ & MAE   \\ \hline 		
			$S^{D}$&$\beta=6 \cdot 10^{-2}$ & 0.0079& 0.2361& 3.18 &$3.3670\cdot 10^{-5}$&$3.7239\cdot 10^{-5}$& -  &247.6774  &247.7818 & 253.2188 & $3.7709 \cdot 10^{-6}$& $\boldsymbol{0.0105}$   \\
			outliers & & 0.0070 &0.2362& 3.18  &$3.3138 \cdot 10^{-5}$ &$3.6652\cdot 10^{-5}$ & - & 247.4937 & 247.5963 & 253.0687 & $3.6496 \cdot 10^{-6}$&$\boldsymbol{0.0215} $ \\	\\
			MDLE&$\beta=6 \cdot 10^{-2}$ & 0.0064 & 0.2341 & 3.0193  &$3.9113 \cdot 10^{-5}$&$1.0933  \cdot 10^{-4}$ & 0.0411 &249.4574 &  249.6679 & 257.7694&$4.0746 \cdot 10^{-7}$&0.0117   \\
			outliers & &  0.0055 &  0.2340 & 3.0445  &$3.7408\cdot 10^{-5}$ &$1.0433 \cdot 10^{-4}$ & 0.0408 & 250.1198 & 250.3267  &258.4823 &$3.7776 \cdot 10^{-7}$ &0.0239 \\	  \hline
			EEs, $\rho$ & TC & $\hat{\mu}$ & $\hat{\sigma}$  &$\alpha, \hat{\alpha}$ & $\text{Var}_{S_q}(\hat{\mu})$ & $\text{Var}_{S_q}(\hat{\sigma})$ & $\text{Var}_{S_q}(\hat{\alpha})$ & $\text{AIC}_{S_q}$ & $\text{cAIC}_{S_q}$ &  $\text{BIC}_{S_q}$ & $\text{Vol}_q(S_q;\hat{\mu},\hat{\sigma}, \hat{\alpha})$ & MAE   \\ \hline 
			$S_q$& $q=.93$ & 0.0095& 0.2550 & 3.18 &$3.5113\cdot 10^{-5}$ &3$.8234\cdot 10^{-5}$ & -  & 259.2036 & 259.3080 & 264.7450  & $3.9020 \cdot 10^{-6}$&0.0131   \\
			outliers &  & 0.0096 &0.2579& 3.18  &$3.5357\cdot 10^{-5}$ &$3.8487\cdot 10^{-5}$ &- &255.6675 & 255.7701  &261.2425 & $3.8630 \cdot 10^{-6}$&0.0225 \\ \\	 
			MqLE& $q=.93$ & 0.0072&0.2472 &2.9071 &$4.4251\cdot 10^{-5}$ &$1.2744 \cdot 10^{-4}$& 0.0314 & 265.4870 &  265.6976 & 273.7991&$4.0498 \cdot 10^{-7}$&0.0121 \\
			outliers &  & 0.0074 &0.2548 & 2.9019 &$4.6635 \cdot 10^{-5}$ &$1.3427\cdot 10^{-4}$ & 0.0309 & 257.2779 & 257.4848 & 265.6404&$4.1279 \cdot 10^{-7}$&0.0234 \\	\hline
			
		\end{tabular}
	}
\end{table}

%% file: table5.tex
\begin{table}[!htb]
	\centering
	\caption{$1^{st}$ design}
	\label{sim1Design}
	\scalebox{0.65}{
		\begin{tabular}{c|c|cccc|cccc|cccc}\hline
			EEs&$\boldsymbol \theta$ & TC$_1$ & $\hat{\boldsymbol \theta}=(\hat{\mu},\hat{\sigma})$ & $\widehat{\text{Var}}(\hat{\boldsymbol \theta})$ & $\widehat{\text{MSE}}(\hat{\boldsymbol \theta})$  & TC$_2$ & $\hat{\boldsymbol \theta}$ & $\widehat{\text{Var}}(\hat{\boldsymbol \theta})$ & $\widehat{\text{MSE}}(\hat{\boldsymbol \theta})$ & TC$_3$ &$\hat{\boldsymbol \theta}$  & $\widehat{\text{Var}}(\hat{\boldsymbol \theta})$ &$\widehat{\text{MSE}}(\hat{\boldsymbol \theta})$ \\  \hline
			&  &$-k,t$ & \multicolumn{3}{c}{$n=110$}  & $-k,t$   &\multicolumn{3}{c}{$n=210$}  &$-k,t$ &  \multicolumn{3}{c}{$n=410$}  \\ \hline 			
			\multirow{2}{*}{$S^H$}& & \multirow{2}{*}{-1.38,1.38}& 0.0852 & 0.0025 & 0.0097 & \multirow{2}{*}{-2.47,2.47} & 0.0718 & 0.0014  & 0.0065 & \multirow{2}{*}{-5.01,5.01} & 0.0638 & 0.0008 &  0.0048   \\
			&   \multirow{8}{*}{$\mu=0$} &  &  1.0058 &   0.0080 &  \it{0.0080}& & 1.0089  &  0.0059 &   0.0060 & & 1.0001   & 0.0050 &   0.0050  \\
			\multirow{2}{*}{$S^{*H}$} &\multirow{9}{*}{$\sigma=1$}  &\multirow{2}{*}{-0.79,0.78} & 0.0830  & 0.0019 & \it{0.0088}  & \multirow{2}{*}{-1.02,1.01}& 0.0438 & 0.0009& 0.0028& \multirow{2}{*}{-1.25,1.24} & 0.0236 &  0.0004 &   0.0010     \\ 
			& &    &  1.0046  &  0.0084 &   0.0084   & & 1.0041   & 0.0044  &  0.0044&   &  1.0054 &   0.0020  &  0.0020 \\
			\multirow{2}{*}{$S^{*}$}&   & \multirow{2}{*}{-0.37,0.36} &  0.1016  &  0.0019&    0.0122   & \multirow{2}{*}{-1.04,1.03} & \it{0.0435}  &  0.0009  &  \it{0.0028} &   \multirow{2}{*}{-1.47,1.46} & 0.0167  &  0.0005&   \it{0.0007}   \\
			& &    &  1.1625 &   0.0107 &   0.0371   & & \it{1.0035}  &  0.0044  &  \it{0.0044} &   &  1.0058 &   0.0019   & \it{0.0019}  \\ \\
			\multirow{2}{*}{$S_q$}	   & & 	\multirow{2}{*}{$q=0.84 $}  & 0.0190  & 0.0019 &  0.0023  &  \multirow{2}{*}{0.875} & 0.0142 &  0.0009  &  0.0011 & \multirow{2}{*}{$q=0.905$} &0.0108 & 0.0004  &0.0006  \\
			&&& 1.0107  &  0.0102  &  0.0103  & & 1.0066 &   0.0040   & 0.0041 &  &  1.0075 &  0.0017 &    0.0018 \\
			\multirow{2}{*}{$S^{D}$} & & \multirow{2}{*}{$\beta=3\cdot10^{-3} $}  &  0.0062 &  0.0020 & \bf{0.0021}  & \multirow{2}{*}{$\beta=10^{-3} $} & 0.0043  & 0.0010& \bf{0.0010} &  \multirow{2}{*}{$\beta=0.8\cdot10^{-3}$} & 0.0021 & 0.0005 & \bf{0.0005}   \\
			& & &  1.0050   & 0.0061&    \bf{0.0061}&     & 1.0079 &   0.0028 & \bf{0.0029}   & & 1.0010  &  0.0014 &   \bf{0.0014}  \\ 			\hline
			$\rho$&$\boldsymbol \theta$ & TC$_1$ & $\hat{\boldsymbol \theta}=(\hat{\mu},\hat{\sigma},\hat{\alpha})$ & $\widehat{\text{Var}}(\hat{\boldsymbol \theta})$ & $\widehat{\text{MSE}}(\hat{\boldsymbol \theta})$  & TC$_2$ & $\hat{\boldsymbol \theta}$ & $\widehat{\text{Var}}(\hat{\boldsymbol \theta})$ & $\widehat{\text{MSE}}(\hat{\boldsymbol \theta})$ & TC$_3$ &$\hat{\boldsymbol \theta}$  & $\widehat{\text{Var}}(\hat{\boldsymbol \theta})$ &$\widehat{\text{MSE}}(\hat{\boldsymbol \theta})$ \\  \hline
			\multirow{3}{*}{MqLE}&\multirow{4}{*}{ $\mu=0$}&\multirow{3}{*}{$q=0.625$}&0.0069&0.0016&\bf{0.0017} &\multirow{3}{*}{$q=0.625$} & 0.0030   & 0.0006 &  \bf{0.0006} &\multirow{3}{*}{$q=0.6$}  & 0.0016   & 0.0003&    \bf{0.0003} \\
			&\multirow{4}{*}{  $\sigma=1$} & &0.7589 & 0.0325 & 0.0906& & 0.7766 &   0.0142 &   0.0641 & & 0.7630  &  0.0075&    0.0637  \\
			&	\multirow{4}{*}{ $\alpha_2=2$} & & 2.0080 &   0.4568 & 0.4568 &   &  2.0013 &   0.1913 &    0.1913 &  & 1.9871   & 0.1036 &   0.1037 \\
			\multirow{3}{*}{MDLE}& &\multirow{3}{*}{$\beta=2.5\cdot10^{-3}$} & 0.0098 &   0.0021  &  0.0022 &\multirow{3}{*}{$\beta=1.7 \cdot 10^{-3}$} &0.0048 &   0.0010 &   0.0010& \multirow{3}{*}{$\beta=1.8 \cdot 10^{-3} $} &   0.0026  &  0.0005   & 0.0005 \\
			&&& 0.9734  & 0.0261& \bf{0.0268} & & 0.9936  &  0.0098  & \bf{0.0099}  & & 0.9979  &  0.0052  &  \bf{0.0052} \\
			&&& 2.0282  & 0.3786& \bf{0.3794}  & &2.0402 &   0.1389 &  \bf{0.1405}& &     2.0597 &   0.0709&   \bf{0.0745} \\ 	\hline
		\end{tabular}
	}
\end{table}

%% file: table6.tex
\begin{table}[!htb]
	\centering
	\caption{$2^{nd}$ design}
	\label{sim2Design}
	\scalebox{0.65}{
		\begin{tabular}{c|c|cccc|cccc|cccc}\hline
			EEs&$\boldsymbol \theta$ & TC$_1$ & $\hat{\boldsymbol \theta}=(\hat{\mu},\hat{\sigma})$ & $\widehat{\text{Var}}(\hat{\boldsymbol \theta})$ & $\widehat{\text{MSE}}(\hat{\boldsymbol \theta})$  & TC$_2$ & $\hat{\boldsymbol \theta}$ & $\widehat{\text{Var}}(\hat{\boldsymbol \theta})$ & $\widehat{\text{MSE}}(\hat{\boldsymbol \theta})$ & TC$_3$ &$\hat{\boldsymbol \theta}$  & $\widehat{\text{Var}}(\hat{\boldsymbol \theta})$ &$\widehat{\text{MSE}}(\hat{\boldsymbol \theta})$ \\  \hline
			&  &$-k,t$ & \multicolumn{3}{c}{$n=110$}  & $-k,t$   &\multicolumn{3}{c}{$n=210$}  &$-k,t$ &  \multicolumn{3}{c}{$n=410$}  \\ \hline 			
			\multirow{2}{*}{$S^H$}& & \multirow{2}{*}{-1.72,1.72}&   0.0803  &  0.0018  &  0.0083 & \multirow{2}{*}{-3.31,3.31} &     0.0723  &  0.0009 &  0.0061 & \multirow{2}{*}{-7.85,7.85} & 0.0662   & 0.0006  &  0.0050  \\
			& \multirow{8}{*}{$\mu=0$}  &  &  1.0062  &  0.0123 &  \it{0.0123} & & 1.0025 & 0.0117&0.0117   & & 1.0056  &  0.0158  & 0.0158  \\
			\multirow{2}{*}{$S^{*H}$} & \multirow{9}{*}{$\sigma=1$} &\multirow{2}{*}{-0.84,0.84}&0.0701   & 0.0023 &   \it{0.0072}  & \multirow{2}{*}{-1.05,1.05}& 0.0312& 0.0012& \it{0.0022}& \multirow{2}{*}{-1.22,1.22} &   0.0142  &  0.0005  &  0.0007    \\
			& &    &        0.9986  &  0.0155 &   0.0155   & & 1.0033  &  0.0064 & 0.0064 &   &     1.0071 &   0.0022   & 0.0022 \\
			\multirow{2}{*}{$S^{*}$}&   & \multirow{2}{*}{-0.57,0.58} &  0.0868 &   0.0019& 0.0094 & \multirow{2}{*}{-1.07,1.09} &     0.0330 &   0.0012  &  0.0023 &   \multirow{2}{*}{-1.28,1.29} &  0.0121   & 0.0006 &    \it{0.0007}  \\
			& &    &  1.0629  &  0.0173 &   0.0213   & & 1.0025   & 0.0060  &  \it{0.0060} &   &  1.0040  &  0.0018   & \it{0.0018}\\ \\
			\multirow{2}{*}{$S_q$}	  & & 	\multirow{2}{*}{$q=0.86 $}  &0.0144  &  0.0022   &0.0024   &  \multirow{2}{*}{0.9} & 0.0083  &  0.0011  &  0.0011 & \multirow{2}{*}{$q=0.93$} &    0.0047   & 0.0005 &   0.0006 \\
			&&& 1.0028  &  0.0042 &   0.0042  & & 1.0009 &   0.0021 &   0.0021 &  &     1.0034  &  0.0011  &  0.0011 \\
			\multirow{2}{*}{$S^{D}$} & & \multirow{2}{*}{$\beta=7\cdot{10^{-3}} $}  &  0.0139   & 0.0023  &  \bf{0.0024} & \multirow{2}{*}{$\beta=4\cdot10^{-3} $} &     0.0077  &  0.0011 &   \bf{0.0011} &  \multirow{2}{*}{$\beta=1.5\cdot10^{-3}$} & 0.0044  &  0.0005 &  \bf{0.0006}    \\
			& & &     1.0068 &   0.0042 &  \bf{0.0042} &     & 1.0017  &  0.0019 &  \bf{0.0019} & &   1.0040  &  0.0009  & \bf{0.0010}  \\ 			\hline
			$\rho$&$\boldsymbol \theta$ & TC$_1$ & $\hat{\boldsymbol \theta}=(\hat{\mu},\hat{\sigma},\hat{\alpha})$ & $\widehat{\text{Var}}(\hat{\boldsymbol \theta})$ & $\widehat{\text{MSE}}(\hat{\boldsymbol \theta})$  & TC$_2$ & $\hat{\boldsymbol \theta}$ & $\widehat{\text{Var}}(\hat{\boldsymbol \theta})$ & $\widehat{\text{MSE}}(\hat{\boldsymbol \theta})$ & TC$_3$ &$\hat{\boldsymbol \theta}$  & $\widehat{\text{Var}}(\hat{\boldsymbol \theta})$ &$\widehat{\text{MSE}}(\hat{\boldsymbol \theta})$ \\  \hline
			\multirow{3}{*}{MqLE}&\multirow{4}{*}{ $\mu=0$}&\multirow{3}{*}{$q=0.55$}&    0.0106 &   0.0019  & \bf{0.0020} &\multirow{3}{*}{$q=0.52$}   &   0.0053 &   0.0009 &  \bf{0.0009} &\multirow{3}{*}{$q=0.6$}  & 0.0032  &  0.0004  &  \bf{0.0004} \\
			&\multirow{4}{*}{  $\sigma=1$} & &0.7703&    0.0209&  0.0736& &0.7632  &  0.0102&    0.0663 &  &   0.8111 &   0.0043 &   0.0400\\
			&\multirow{4}{*}{ $\alpha_2=3$} & &  2.7923 & 0.9898 &   1.0329 & &   2.7638 &   0.5645 &   0.6203 & &     2.7337&    0.2323  &  0.3033  \\
			\multirow{3}{*}{MDLE}& &\multirow{3}{*}{$\beta=9\cdot10^{-3}$} & 0.0139  &  0.0021  &0.0022 &\multirow{3}{*}{$\beta=9.5\cdot10^{-3}$}  &0.0074 &   0.0010& 0.0010 & \multirow{3}{*}{$\beta=9\cdot10^{-3}$} &     0.0041 &   0.0005  &  0.0005 \\
			&&& 0.9751 &   0.0087  &  \bf{0.0093} & & 0.9704 &  0.0046 &  \bf{0.0055} & &     0.9711   & 0.0025 &  \bf{0.0033}\\
			&&& 3.0107 &   0.7184 & \bf{0.7185}& &  3.0169&    0.3717  &  \bf{0.3719} & &     3.0141  &  0.2003  & \bf{0.2005} \\ 	\hline
		\end{tabular}
	}
\end{table}

%% file: table7.tex
\begin{table}[!htb]
	\centering
	\caption{$3^{rd}$ design}
	\label{sim3Design}
	\scalebox{0.65}{
		\begin{tabular}{c|c|cccc|cccc|cccc}\hline
			EEs&$\boldsymbol \theta$ & TC$_1$ & $\hat{\boldsymbol \theta}=(\hat{\mu},\hat{\sigma})$ & $\widehat{\text{Var}}(\hat{\boldsymbol \theta})$ & $\widehat{\text{MSE}}(\hat{\boldsymbol \theta})$  & TC$_2$ & $\hat{\boldsymbol \theta}$ & $\widehat{\text{Var}}(\hat{\boldsymbol \theta})$ & $\widehat{\text{MSE}}(\hat{\boldsymbol \theta})$ & TC$_3$ &$\hat{\boldsymbol \theta}$  & $\widehat{\text{Var}}(\hat{\boldsymbol \theta})$ &$\widehat{\text{MSE}}(\hat{\boldsymbol \theta})$ \\  \hline
			&  &$-k,t$ & \multicolumn{3}{c}{$n=110$}  & $-k,t$   &\multicolumn{3}{c}{$n=210$}  &$-k,t$ &  \multicolumn{3}{c}{$n=410$}  \\ \hline 			
			\multirow{2}{*}{$S^H$}& & \multirow{2}{*}{-1.34,1.34}& 0.0773  &  0.0018  &  0.0078 & \multirow{2}{*}{-2.53,2.53} &0.0676  &  0.0010 &  0.0056 & \multirow{2}{*}{-5.21,5.21} &     0.0588 &   0.0005 &   0.0040   \\
			& \multirow{8}{*}{$\mu=0$}  &  & 1.0005   & 0.0219 &  0.0219& &  1.0022  &  0.0186 & 0.0186  & &    0.9984 &   0.0200 &   0.0200 \\
			\multirow{2}{*}{$S^{*H}$} & \multirow{9}{*}{$\sigma=1$}  &\multirow{2}{*}{-1.13,1.12} & \it{0.0012} & 0.0024& \it{0.0024}  & \multirow{2}{*}{-1.29,1.28}& 0.0020 &   0.0012&    0.0012& \multirow{2}{*}{-1.40,1.38} &     0.0004  &  0.0005  &  \it{0.0005}    \\
			&  &    &      1.0027 &   0.0105 & 0.0105   & & 1.0107   & 0.0035 &   0.0036&   &     1.0019 &   0.0013  &  0.0013 \\
			\multirow{2}{*}{$S^{*}$}&  & \multirow{2}{*}{-1.2,1.19} & \it{0.0042} &   0.0024  &  \it{0.0024}   & \multirow{2}{*}{-1.37,1.35} &  \bf{0.0020}  &  0.0011   & \bf{\it{0.0011}} &   \multirow{2}{*}{-1.48,1.47} &       0.0018  &  0.0006  &  0.0006  \\
			& &    &  1.0074  &  0.0091  &  \it{0.0091}  & &     1.0053   & 0.0029 &  \it{0.0030} &   &      1.0019  &  0.0012 &   \it{0.0012}\\ \\
			\multirow{2}{*}{$S_q$}  & &\multirow{2}{*}{$q=0.85 $}  & 0.0184  &  0.0023  & 0.0026   &  \multirow{2}{*}{$q=0.89$} &0.0122  &  0.0011  &  0.0012 & \multirow{2}{*}{$q=0.92$} &    0.0084 &   0.0005 &   0.0006  \\
			&&&     1.0053 &   0.0051 &   0.0051  & &     1.0046 &   0.0024 &   0.0024 &  &      1.0045  &  0.0011 &   0.0011 \\
			\multirow{2}{*}{$S^{D}$} & & \multirow{2}{*}{$\beta=4\cdot{10^{-3}} $}  &  0.0096  &  0.0023  &  \bf{0.0023} & \multirow{2}{*}{$\beta=2\cdot10^{-3} $} &   \bf{0.0051} &   0.0011 &  \bf{0.0011} &  \multirow{2}{*}{$\beta=10^{-3}$} &    0.0028   & 0.0005  &  \bf{0.0005}    \\
			& & & 1.0076 &   0.0041  & \bf{0.0042}&     &     1.0038 &   0.0020 &  \bf{0.0020} & &      1.0026 &   0.0009 &  \bf{0.0009}  \\ 			\hline
			$\rho$&$\boldsymbol \theta$ & TC$_1$ & $\hat{\boldsymbol \theta}=(\hat{\mu},\hat{\sigma},\hat{\alpha})$ & $\widehat{\text{Var}}(\hat{\boldsymbol \theta})$ & $\widehat{\text{MSE}}(\hat{\boldsymbol \theta})$  & TC$_2$ & $\hat{\boldsymbol \theta}$ & $\widehat{\text{Var}}(\hat{\boldsymbol \theta})$ & $\widehat{\text{MSE}}(\hat{\boldsymbol \theta})$ & TC$_3$ &$\hat{\boldsymbol \theta}$  & $\widehat{\text{Var}}(\hat{\boldsymbol \theta})$ &$\widehat{\text{MSE}}(\hat{\boldsymbol \theta})$ \\  \hline
			\multirow{3}{*}{MqLE}&\multirow{4}{*}{ $\mu=0$}&\multirow{3}{*}{$q=0.65$}&  0.0080 &   0.0015  & \bf{0.0016} &\multirow{3}{*}{$q=0.65$} & 0.0042  &  0.0007 &  \bf{0.0008} & \multirow{3}{*}{$q=0.65$}  &   0.0024 &   0.0004 &   \bf{0.0004} \\
			&\multirow{4}{*}{  $\sigma=1$} & &    0.8196 &   0.0162 &   0.0488& &    0.8284 &   0.0076  &  0.0370 & & 0.8317 &   0.0042 &   0.0325 \\
			&	\multirow{4}{*}{ $\alpha_2=3$} & &   2.7331  &  0.7868  &  0.8581 & &     2.7371 &   0.3979  &  0.4670 & & 2.7300  &  0.2200  &  0.2929 \\
			\multirow{3}{*}{MDLE}& &\multirow{3}{*}{$\beta=7\cdot10^{-3}$} &0.0100 &   0.0020  &  0.0021 & \multirow{3}{*}{$\beta=7\cdot10^{-3}$} &   0.0050  &  0.0010 &   0.0010 & \multirow{3}{*}{$\beta=8\cdot10^{-3}$} &     0.0028  &  0.0005 &   0.0005\\
			&&& 0.9701&    0.0087   & \bf{0.0096}& &0.9700    &0.0046  &  \bf{0.0055} & &     0.9690  &  0.0026 &  \bf{0.0036}  \\
			&&&3.0018   & 0.7229   & \bf{0.7229} &  &   2.9949  &  0.3575 &  \bf{0.3576} & &     3.0022 &   0.2049  &  \bf{0.2049} \\ 	\hline
		\end{tabular}
	}
\end{table}

%% file: table8.tex
\begin{table}[!htb]
	\centering
	\caption{$4^{th}$ design}
	\label{sim4Design}
	\scalebox{0.65}{
		\begin{tabular}{c|c|cccc|cccc|cccc}\hline
			EEs&$\boldsymbol \theta$ & TC$_1$ & $\hat{\boldsymbol \theta}=(\hat{\mu},\hat{\sigma})$ & $\widehat{\text{Var}}(\hat{\boldsymbol \theta})$ & $\widehat{\text{MSE}}(\hat{\boldsymbol \theta})$  & TC$_2$ & $\hat{\boldsymbol \theta}$ & $\widehat{\text{Var}}(\hat{\boldsymbol \theta})$ & $\widehat{\text{MSE}}(\hat{\boldsymbol \theta})$ & TC$_3$ &$\hat{\boldsymbol \theta}$  & $\widehat{\text{Var}}(\hat{\boldsymbol \theta})$ &$\widehat{\text{MSE}}(\hat{\boldsymbol \theta})$ \\  \hline
			&  &$-k,t$ & \multicolumn{3}{c}{$n=110$}  & $-k,t$   &\multicolumn{3}{c}{$n=210$}  &$-k,t$ &  \multicolumn{3}{c}{$n=410$}  \\ \hline 			
			\multirow{2}{*}{$S^H$}& & \multirow{2}{*}{-1.07,1.07}& 0.0931   & 0.0031 &   0.0118 & \multirow{2}{*}{-1.38,1.38} & 0.0532  &  0.0017  &  0.0045 & \multirow{2}{*}{-1.73,1.73} & 0.0308  &  0.0009 &   0.0019  \\
			& \multirow{8}{*}{$\mu=0$}  &  &     1.0083    &0.0120 &   0.0120& &    1.0029  &  0.0055 &   0.0056 & &    1.0011&    0.0027   & 0.0027  \\
			\multirow{2}{*}{$S^{*H}$} &\multirow{9}{*}{$\sigma=1$}&\multirow{2}{*}{-1.21,1.22} & 0.1001 &   0.0016 &   0.0116  & \multirow{2}{*}{-1.47,1.48}&  0.0620  &  0.0005 &   0.0044& \multirow{2}{*}{-1.74,1.76} &    0.0397  &  0.0002 &   0.0018     \\
			&  &    &  1.0004   & 0.0101 &   \it{0.0101}   & &    1.0016   & 0.0050 &   \it{0.0050}&   &       1.0043  &  0.0024 &   0.0024 \\
			\multirow{2}{*}{$S^{*}$}&   & \multirow{2}{*}{-1.76,1.77} & 0.0685  &  0.0014  &  \it{0.0060}   & \multirow{2}{*}{-2.34,2.36} &0.0331   & 0.0005 &   \it{0.0016} &   \multirow{2}{*}{-2.71,2.73} &     0.0175   & 0.0002 &   \it{0.0005}   \\
			& &    &     1.0089  &  0.0105 &   0.0106   & & 1.0098  &  0.0052  &  0.0053 &   &     0.9964  &  0.0024  &  \it{0.0024}  \\ \\
			\multirow{2}{*}{$S_q$}	   & & 	\multirow{2}{*}{$q=0.815$}  &0.0502 &   0.0011 &   0.0036 &  \multirow{2}{*}{$q=0.875$} & 0.0290 &   0.0004   & 0.0012 & \multirow{2}{*}{$q=0.918$} &    0.0165 &   0.0001  &  0.0004 \\
			&&& 1.0046  &  0.0088  &  \bf{0.0089}  & &     1.0043 &    0.0040 &   \bf{0.0040} &  &      1.0029  &  0.0019 &   \bf{0.0019}\\
			\multirow{2}{*}{$S^{D}$} & & \multirow{2}{*}{$\beta=9\cdot{10^{-3}} $}  &     0.0400  &  0.0013 &  \bf{0.0029} & \multirow{2}{*}{$\beta=4\cdot10^{-3} $} &     0.0231    &0.0004 &  \bf{0.0010} &  \multirow{2}{*}{$\beta=2\cdot{10^{-3}}$} &    0.0129   & 0.0002  &  \bf{0.0003}  \\
			& & & 1.0041 &   0.0104  &  0.0104&     &     1.0070 &   0.0046 &   0.0047 & &      1.0018   & 0.0021  &  0.0021  \\ 			\hline
			$\rho$&$\boldsymbol \theta$ & TC$_1$ & $\hat{\boldsymbol \theta}=(\hat{\mu},\hat{\sigma},\hat{\alpha})$ & $\widehat{\text{Var}}(\hat{\boldsymbol \theta})$ & $\widehat{\text{MSE}}(\hat{\boldsymbol \theta})$  & TC$_2$ & $\hat{\boldsymbol \theta}$ & $\widehat{\text{Var}}(\hat{\boldsymbol \theta})$ & $\widehat{\text{MSE}}(\hat{\boldsymbol \theta})$ & TC$_3$ &$\hat{\boldsymbol \theta}$  & $\widehat{\text{Var}}(\hat{\boldsymbol \theta})$ &$\widehat{\text{MSE}}(\hat{\boldsymbol \theta})$ \\  \hline
			\multirow{3}{*}{MqLE}&\multirow{4}{*}{ $\mu=0$}&\multirow{3}{*}{$q=0.525$}&  0.0180 &0.0042&\bf{0.0045}&\multirow{3}{*}{$q=0.565$} &  0.0085 &   0.0011  & \bf{0.0012}&\multirow{3}{*}{$q=0.605$} &  0.0045   & 0.0003  & \bf{0.0003} \\
			&\multirow{4}{*}{  $\sigma=1$} & & 0.5755&0.0617&0.2419& &     0.6351  &  0.0320   & 0.1652 & &     0.6824 &   0.0161  &  0.1170\\
			&	\multirow{4}{*}{ $\alpha_2=1.3$} & &1.2940 &   0.2857  &  0.2857 & &     1.2907 &   0.0937 &   0.0938 & & 1.3030  &  0.0424  &  0.0424\\
			\multirow{3}{*}{MDLE}& &\multirow{3}{*}{$\beta=6\cdot10^{-3}$} &     0.0488  &  0.0022  &  0.0045 &\multirow{3}{*}{$\beta=3\cdot10^{-3}$}&  0.0275 &   0.0007 &   0.0014 & \multirow{3}{*}{$\beta=1.2\cdot10^{-3}$}&     0.0161 &   0.0003    &0.0005 \\
			&&&     1.0139 &   0.0548  &  \bf{0.0550} & &    1.0137 &   0.0272 &  \bf{0.0273} & &     1.0146 &   0.0150  & \bf{0.0152}\\
			&&&     1.3089 &   0.1204  &  \bf{0.1205} & &     1.3022  &  0.0449 &  \bf{0.0449} & &     1.3009    &0.0233    &\bf{0.0233} \\ 	\hline
		\end{tabular}
	}
\end{table}

%% file: EEq.bbl
\begin{thebibliography}{99}
	
\bibitem{PardoSD} L. Pardo, Statistical inference based on divergence measures, CRC Press, Taylor \& Francis Group, 2005.

\bibitem{Tsallisbook09} C. Tsallis, Introduction to Nonextensive Statistical Mechanics: Approaching a Complex World, Springer, New York, 2009.

\bibitem{Tsallis88} Tsallis, C. (1988). Possible generalization of Boltzmann-Gibbs statistics. Journal of statistical physics, 52(1-2), 479-487.

\bibitem{AbeOka01book}  Abe, S.,  Okamoto, Y. (Eds.). (2001). Nonextensive statistical mechanics and its applications (Vol. 560). Springer Science and Business Media.

\bibitem{LehmannCas98}  E.L. Lehmann, G. Casella, Theory of point estimation,  Wadsworth \& Brooks/Cole. Pacific Grove, CA, 589, USA, 1998.


\bibitem{FerrariYang10}  Ferrari, D.,  Yang, Y. (2010). Maximum Lq-likelihood estimation. The Annals of Statistics, 38(2), 753-783.


\bibitem{CanKor18}  \c{C}ankaya, M. N.,  Korbel, J. (2018). Least informative distributions in maximum q-log-likelihood estimation. Physica A: Statistical Mechanics and its Applications, 509, 140-150.


\bibitem{MLqEBio} Hasegawa, Y.,  Arita, M. (2009). Properties of the maximum q-likelihood estimator for independent random variables. Physica A: Statistical Mechanics and its Applications, 388(17), 3399-3412.


\bibitem{FerrariOpere} Giuzio, M., Ferrari, D.,  Paterlini, S. (2016). Sparse and robust normal and t-portfolios by penalized Lq-likelihood minimization. European Journal of Operational Research, 250(1), 251-261.


\bibitem{GentEntReview1} Hanel, R.,  Thurner, S. (2011). A comprehensive classification of complex statistical systems and an axiomatic derivation of their entropy and distribution functions. EPL (Europhysics Letters), 93(2), 20006.


\bibitem{JizbaArimitsu}   Jizba, P.,  Arimitsu, T. (2004). Observability of R\'{e}nyi's entropy. Physical Review E, 69(2), 026128.

\bibitem{Bercher12a} Bercher, J. F. (2012). A simple probabilistic construction yielding generalized entropies and divergences, escort distributions and q-Gaussians. Physica A: Statistical Mechanics and its Applications, 391(19), 4460-4469.

\bibitem{Jizba3}  P. Jizba, Information theory and generalized statistics, in: H.-T.~Elze, ed., Decoherence and Entropy in Complex Systems, Lecture Notes in Physics, vol. 633, Springer-Verlag, Berlin, 2003, p.362.


\bibitem{JizKorhybrid}  Jizba, P.,  Korbel, J. (2016). On q-non-extensive statistics with non-Tsallisian entropy. Physica A: Statistical Mechanics and its Applications, 444, 808-827.


\bibitem{Basuetal98}  Basu, A., Harris, I. R., Hjort, N. L.,  Jones, M. C. (1998). Robust and efficient estimation by minimising a density power divergence. Biometrika, 85(3), 549-559.


\bibitem{Vajda86} Vajda, I. (1986). Efficiency and robustness control via distorted maximum likelihood estimation. Kybernetika, 22(1), 47-67.


\bibitem{Hub64}	Huber, P. J. (1964). Robust estimation of a location parameter: Annals Mathematics Statistics, 35.

\bibitem{God60}  Godambe, V. P. (1960). An optimum property of regular maximum likelihood estimation. The Annals of Mathematical Statistics, 31(4), 1208-1211.


\bibitem{Omer98} \"{O}zt\"{u}rk, \"{O}. (1998). Theory and Methods: A Robust and Almost Fully Efficient M‐Estimator. Australian and New Zealand Journal of Statistics, 40(4), 415-424.

\bibitem{Hampeletal86} F.R. Hampel, E. M. Ronchetti, P. J. Rousseeuw, W. A. Stahel, Robust statistics: The approach based on influence functions, Wiley Series in Probability and Statistics, New York, 1986.


\bibitem{SmithODE}  Smith, H. L. (2011). An introduction to delay differential equations with applications to the life sciences (Vol. 57). New York: Springer.

\bibitem{infocomplestamod} J. Rissanen, Information and Complexity in Statistical Modeling, Springer, New York, 2007. 
\bibitem{MyungBalasubramanian}   Myung, I. J., Balasubramanian, V.,  Pitt, M. A. (2000). Counting probability distributions: Differential geometry and model selection. Proceedings of the National Academy of Sciences, 97(21), 11170-11175.


\bibitem{ClusteringFisher} Gattone, S. A., De Sanctis, A., Russo, T.,  Pulcini, D. (2017). A shape distance based on the Fisher-Rao metric and its application for shapes clustering. Physica A: Statistical Mechanics and its Applications, 487, 93-102. 


\bibitem{Maybankvol}  Maybank, S. J. (2006). Application of the Fisher-Rao metric to structure detection. Journal of Mathematical Imaging and Vision, 25(1), 49-62.

\bibitem{AmariIG}  Amari, S., Information Geometry and Its Applications, in: Applied Mathematical Sciences, Springer, 2016.

\bibitem{Shannon}  Shannon, Claude Elwood. "A mathematical theory of communication." ACM SIGMOBILE mobile computing and communications review 5.1 (2001): 3-55.


\bibitem{Ernst12} Ernst, T. (2012). A comprehensive treatment of q-calculus. Springer Science \& Business Media.

\bibitem{Abe97} Abe, S. (1997). A note on the q-deformation-theoretic aspect of the generalized entropies in nonextensive physics. Physics Letters A, 224(6), 326-330.

\bibitem{Jackson08} F. H. Jackson, On q-functions and a certain difference operator, Appl. Math. Sci, vol. 46, no. 2, 	pp. 253–281, 1908, doi: 10.1017/S0080456800002751.

\bibitem{Wadatwopara}  Wada, T., Suyari, H. (2007). A two-parameter generalization of Shannon–Khinchin axioms and the uniqueness theorem. Physics Letters A, 368(3-4), 199-205.

\bibitem{Beck97}  C. Beck, F. Schl\"{o}gl, Thermodynamics of Chaotic Systems: An Introduction, Cambridge Un. Press, Cambridge, 1997.

\bibitem{Beck04} Beck, C. (2004). Superstatistics, escort distributions, and applications. Physica A: Statistical Mechanics and its Applications, 342(1-2), 139-144.

\bibitem{Harte01} D. Harte, Multifractals, Theory and Applications, Chapman \& Hall/CRC, London, 2001.	



\bibitem{KullbackIT} S. Kullback, Information Theory and Statistics, Courier Corporation, USA, 1997.	


\bibitem{Plastinoetal97}     Plastino, A., Plastino, A. R.,  Miller, H. G. (1997). Tsallis nonextensive thermostatistics and Fisher's information measure. Physica A: Statistical Mechanics and its Applications, 235(3-4), 577-588.	


\bibitem{CanEnt18} \c{C}ankaya, M. N. (2018). Asymmetric bimodal exponential power distribution on the real line. Entropy, 20(1), 23.

\bibitem{Cantez19} \c{C}ankaya, M. N.,  Arslan, O. (2020). On the robustness properties for maximum likelihood estimators of parameters in exponential power and generalized T distributions. Communications in Statistics-Theory and Methods, 49(3), 607-630.


\bibitem{Maronna76}  Maronna, R. A. (1976). Robust M-estimators of multivariate location and scatter. The annals of statistics, 51-67.


\bibitem{Bercher12}  Bercher, J. F. (2012). On generalized Cramér–Rao inequalities, generalized Fisher information and characterizations of generalized q-Gaussian distributions. Journal of Physics A: Mathematical and Theoretical, 45(25), 255303.

\bibitem{PhDthesis} \c{C}ankaya, M. N. 2015. M-Estimators with asymmetric influence function: Properties and their applications. PhD diss., University of Ankara.

\bibitem{GodTh84} 	Godambe, V. P., Thompson, M. E. (1984). Robust estimation through estimating equations. Biometrika, 71(1), 115-125.
\bibitem{ProfMary} M. Thompson, e-mail communication.

\bibitem{MalikAr92}  Malik SC, Arora S, 1992. Mathematical analysis. New Age International.

\bibitem{Jan17} Korbel, J. (2017). Rescaling the nonadditivity parameter in Tsallis thermostatistics. Physics Letters A, 381(32), 2588-2592.

\bibitem{Jan18} Korbel, J., Hanel, R., Thurner, S. (2018). Classification of complex systems by their sample-space scaling exponents. New Journal of Physics, 20(9), 093007.

\bibitem{Jan19} Jizba, P., Korbel, J. (2019). Maximum entropy principle in statistical inference: Case for non-Shannonian entropies. Physical review letters, 122(12), 120601.


\bibitem{Daroczy}   Dar\'{o}czy, Z. (1970). Generalized information functions. Information and control, 16(1), 36-51.


\bibitem{Fisher25} Fisher, R. A. (1925, July). Theory of statistical estimation. In Mathematical Proceedings of the Cambridge Philosophical Society (Vol. 22, No. 5, pp. 700-725). Cambridge University Press.


\bibitem{CanKor16}  \c{C}ankaya, M. N.,  Korbel, J. (2017). On statistical properties of Jizba-Arimitsu hybrid entropy. Physica A: Statistical Mechanics and its Applications, 475, 1-10.

\bibitem{Hub81}   Huber, P.J. {\it Robust statistics},  1981, John Wiley and Sons, New York, USA.

\bibitem{GodHeyquasi}  Godambe, V. P.,  Heyde, C. C. (2010). Quasi-likelihood and optimal estimation. In Selected works of cc heyde (pp. 386-399). Springer, New York, NY.

\bibitem{Aka73} H. Akaike,  Information theory and an extension of the maximum likelihood principle.  In B. N. Petrov \& B. F. Csaki (Eds.), Second International Symposium on Information Theory,  Academiai Kiado: Budapest, 267-281, 1973.

\bibitem{Hamparsum87}  Bozdogan, H. (1987). Model selection and Akaike's information criterion (AIC): The general theory and its analytical extensions. Psychometrika, 52(3), 345-370.
\bibitem{Roncpt97} Ronchetti, E. (1997). Robustness aspects of model choice. Statistica Sinica, 327-338.

\bibitem{Elsal09} Elsalloukh, H. 2009. Further results on the epsilon-skew exponential power distribution family.  Far East Journal of Theoretical Statistics 28:201–12.

\bibitem{Lucas97} Lucas, A. (1997). Robustness of the student t based M-estimator. Communications in Statistics-Theory and Methods, 26(5), 1165-1182.


\bibitem{Bercher13}  Bercher, J. F. (2013). Some properties of generalized Fisher information in the context of nonextensive thermostatistics. Physica A: Statistical Mechanics and its Applications, 392(15), 3140-3154.



\bibitem{Orkcu15} \"{O}rkc\"{u}, H. H., \"{O}zsoy, V. S., Aksoy, E.,  Do\v{g}an, M. I. (2015). Estimating the parameters of 3-p Weibull distribution using particle swarm optimization: A comprehensive experimental comparison. Applied Mathematics and Computation, 268, 201-226.

\bibitem{temperaturedata}  Available online: legacy.bas.ac.uk:
\begin{verbatim} 
https://legacy.bas.ac.uk/met/READER/surface/Grytviken.All.temperature.txt. 
\end{verbatim}

\bibitem{microarray}  Available online: 
\begin{verbatim} 
https://discover.nci.nih.gov/nature2000/data/selected_data/at_matrix.txt.
\end{verbatim}


\bibitem{matrixproperty} Verhaegen, M.,  Verdult, V. (2007). Filtering and system identification: a least squares approach. Cambridge university press.

\bibitem{Gradshteyn07} I.S. Gradshteyn, I. M. Ryzhik, A. Jeffrey, D. Zwillinger, Table of Integrals, Series, and Products, Sixth Edition, Academic Press, 1171, USA, 2007.
\bibitem{LiYeh12gen}  Li, Y. H.,  Yeh, P. C. (2012). An interpretation of the Moore-Penrose generalized inverse of a singular Fisher Information Matrix. IEEE Transactions on Signal Processing, 60(10), 5532-5536.


\end{thebibliography}
